\documentclass[onecolumn,11pt]{article}

\usepackage[a4paper,margin=1in]{geometry}

\usepackage[affil-it]{authblk}

\usepackage[]{setspace}

\usepackage{rotating}

\makeatletter
\def\@maketitle{%
	\newpage
	\null
	\vskip 2em%
	\begin{center}%
		\let \footnote \thanks
		{\Large\bfseries \@title \par}%
		\vskip 1.5em%
		{\normalsize
			\lineskip .5em%
			\begin{tabular}[t]{c}%
				\@author
			\end{tabular}\par}%
		\vskip 1em%
		{\normalsize \@date}%
	\end{center}%
	\par
	\vskip 1.5em}
\makeatother

\title{The Bicycle Network Improvement Problem}

\author[a,b]{Jisoon Lim}
\author[b]{Kevin Dalmeijer}
\author[c,d]{Subhrajit Guhathakurta}
\author[b]{Pascal Van Hentenryck%
	\thanks{E-mail: \texttt{pascal.vanhentenryck@isye.gatech.edu}; Phone: \texttt{+1 (404) 894-2300}; Address: \texttt{755 Ferst Drive, NW, Atlanta, GA 30332, United States}; Corresponding author}}
\affil[a]{Department of Civil and Environmental Engineering, University of Michigan}
\affil[b]{H. Milton Stewart School of Industrial and Systems Engineering,\protect\\ Georgia Institute of Technology}
\affil[c]{School of City \& Regional Planning, Georgia Institute of Technology}
\affil[d]{Center for Spatial Planning Analytics and Visualization, Georgia Institute of Technology}

\usepackage[latin1]{inputenc} 

\usepackage{mathtools,amsfonts}

\usepackage[english]{babel}

\usepackage{booktabs}
\usepackage{adjustbox}

\usepackage{pdflscape}

\usepackage{enumitem}

\usepackage[short,nocomma]{optidef}

\usepackage[framemethod=tikz]{mdframed}
\newmdenv[innerlinewidth=0.5pt, roundcorner=4pt,linecolor=black,innerleftmargin=6pt,
innerrightmargin=6pt,innertopmargin=6pt,innerbottommargin=6pt]{mybox}

\usepackage{ mathrsfs }

\usepackage{chngpage}

\usepackage[authoryear,round]{natbib}

\usepackage{multirow}

\usepackage{amssymb}
\usepackage{amsmath}
\DeclareMathOperator*{\argmax}{arg\,max}
\DeclareMathOperator*{\argmin}{arg\,min}
\DeclareMathOperator*{\argsup}{arg\,sup}
\DeclareMathOperator*{\arginf}{arg\,inf}

\usepackage{array}
\usepackage{booktabs}
\usepackage{colortbl}

\usepackage{amsthm}

\newtheorem*{myprop*}{Proposition}

\usepackage{float}

\usepackage[hyphens]{url}

\usepackage{xcolor,soul}
\definecolor{dr}{rgb}{0.75,0.00,0.00}
\definecolor{lr}{rgb}{1.00,0.75,0.75}
\sethlcolor{lr}
\newcommand{\fix}[1]{\textbf{\textcolor{dr}{\hl{#1}}}}

\definecolor{dr2}{rgb}{0.00,0.00,0.00}
\definecolor{lr2}{rgb}{0.00,1.00,0.00}
\newcommand{\fixgreen}[1]{\sethlcolor{lr2}\textbf{\textcolor{dr2}{\hl{#1}}}\sethlcolor{lr}} 

\definecolor{dr3}{rgb}{0.00,0.00,1.00}
\newcommand{\changed}[1]{\textbf{\textcolor{dr3}{#1}}}

\usepackage{algorithm}
\usepackage{algpseudocode}
\algdef{SE}[DOWHILE]{Do}{doWhile}{\algorithmicdo}[1]{\algorithmicwhile\ #1}

\usepackage{floatpag}
\floatpagestyle{plain}

\usepackage{afterpage}

\usepackage{tikz}
\usetikzlibrary{decorations.pathmorphing, arrows.meta, shapes, backgrounds}
\usetikzlibrary{arrows, automata}
\usepackage{tkz-euclide}
\usetikzlibrary{calc} 

\usepackage{graphicx}

\usepackage{longtable}

\usepackage{arydshln}

\usepackage{pdflscape}

\usepackage[colorlinks = true,
linkcolor = blue,
urlcolor  = blue,
citecolor = black,
anchorcolor = blue]{hyperref}

\newenvironment{nospacebelowflalign*}
{\setlength{\belowdisplayskip}{0pt}%
	\csname flalign*\endcsname}
{\csname endflalign*\endcsname\ignorespacesafterend}

\newenvironment{nospaceflalign*}
{\setlength{\abovedisplayskip}{0pt}\setlength{\belowdisplayskip}{0pt}%
	\csname flalign*\endcsname}
{\csname endflalign*\endcsname\ignorespacesafterend}

\newcommand{\w}[0]{\omega}
\newcommand{\W}[0]{\Omega}

\newcolumntype{L}[1]{>{\raggedright\let\newline\\\arraybackslash\hspace{0pt}}m{#1}}
\newcolumntype{C}[1]{>{\centering\let\newline\\\arraybackslash\hspace{0pt}}m{#1}}
\newcolumntype{R}[1]{>{\raggedleft\let\newline\\\arraybackslash\hspace{0pt}}m{#1}}

\DeclareMathOperator{\trcl}{trcl}

\begin{document}
	
	\maketitle
	\noindent
	
	\thispagestyle{empty}
	\pagestyle{empty}
	
	\vfill
	
	\begin{spacing}{1.2}
		\begin{abstract}
			\noindent
			Using a bicycle for commuting is still uncommon in US cities, although
			it brings many benefits to both the cyclists and to society as a
			whole.  Cycling has the potential to reduce traffic congestion and
			emissions, increase mobility, and improve public health.  To convince
			people to commute by bike, the infrastructure plays an important role,
			since safety is one of the primary concerns of potential cyclists.
			This paper presents a method to find the best way to improve the
			safety of a bicycle network for a given budget and maximize the number
			of riders that could now choose bicycles for their commuting needs.
			This optimization problem is formalized as the Bicycle Network
			Improvement Problem (BNIP): it selects which roads to improve for a
			set of traveler origin-destination pairs, taking both safety and
			travel distance into account.  The BNIP is modeled as a mixed-integer
			linear program that minimizes a piecewise linear penalty function of
			route deviations of travelers. The MIP is solved using Benders
			decomposition to scale to large instances. The paper also presents an
			in-depth case study for the Midtown area in Atlanta, GA, using actual
			transportation data. The results show that the Benders decomposition
			algorithm allows for solving realistic problem instances and that the
			network improvements may significantly increase the share of bicycles
			as the commuting mode.
			Multiple practical aspects are considered as well, including sequential road improvements, uneven improvement costs, and how to include additional data.
			
			
		\end{abstract}
		\doublespacing
	\end{spacing}
	\thispagestyle{empty}
	\clearpage
	
	\pagestyle{plain}
	\pagenumbering{arabic}
	
	\section{Introduction}
	\label{sec:introduction}
	Using a bicycle for transportation is still uncommon in US cities, but it brings many benefits to both the cyclists and to society as a whole \citep{HandyEtAl2014-PromotingCyclingTransport}.
	Cycling has the potential to reduce traffic congestion and emissions, increase mobility, and improve public health \citep{Northrop2011-BicycleCommuterTrends}.
	Additionally, bikes can serve as an economical alternative to a car, especially for short trips \citep{RyuEtAl2018-TwoStageBicycle}.
	Promoting cycling as an alternative to using a car has been studied extensively, with systematic reviews provided by \citet{OgilvieEtAl2004-PromotingWalkingCycling} and \citet{YangEtAl2010-InterventionsPromoteCycling}.
	
	The benefits of cycling as a mode of transportation have also been recognized by policy makers, and more and more cities have started promoting bicycle usage.
	An example in Atlanta is the \emph{Walk, Bike, Thrive!} plan, which provides a recipe for a more walkable and bikable city \citep{ARC2020-BikePedestrianPlan}.
	Plans like these can play a key role in promoting bicycle usage, as was found by \citet{LanzendorfAnnika2014-CyclingBoomGermany} who studied four German cities.
	Effective cycling policy may also benefit modes that are similar to bicycles, such as e-scooters and e-bikes.
	E-scooters and e-bikes both substitute travel by car \citep{Kroesen2017-WhatExtentDo,Goessling2020-IntegratingEScooters}, and evidence from the Netherlands suggests that car owners are more willing to use e-bikes than conventional bikes \citep{Kroesen2017-WhatExtentDo}.
	
	The low number of cyclists is not due to a lack of interest.
	\citet{DillMcNeil2016-RevisitingFourTypes} questioned 3,000 people in the 50 largest US metropolitan areas about their attitudes towards cycling, and they found that 56\% of the population can be classified as \emph{interested but concerned}.
	One of the key barriers for this group is traffic safety: while most feel comfortable riding on a protected bike lane that is part of a major street, only 16\% would be somewhat comfortable without the bike lane.
	The willingness to cycle is also demonstrated by the surge in US bike ridership during the COVID-19 pandemic \citep{Bryant2020-CyclingexplosionCoronavirus}.
	Many people have started cycling for recreation, but also as a socially-distant alternative to public transit.
	Policy makers hope that this trend continues and that these new cyclists start commuting by bike when they return to work after the pandemic.
	
	To convince people to commute by bike, the infrastructure plays an important role.
	The study by \citet{DillCarr2003-BicycleCommutingFacilities} suggests that, if a city provides the proper infrastructure for cycling, commuters are likely to make use of it.
	\citet{HullOHolleran2014-BicycleInfrastructureCan} study selected European cities to identify whether good design can encourage cycling.
	They found that the design may indeed have a significant impact on mode choice, and that \emph{safety}, \emph{comfort}, and \emph{continuity} were the most influential factors.
	Although safety is especially important to cyclists, safety improvements in the last decades have often focused on motorized vehicles, as highlighted by \citet{CIVITASInitiative2020-SmartChoicesCities} for the European Union.
	This situation can be improved by investing in bicycle infrastructure, which additionally improves the safety for non-cycling road users \citep{Walljasper2016-BicycleCommuterTrends}.
	
	Another important factor, which was not explicitly considered in the previous study, is the \emph{proximity}: the distance between the origin and the destination of the trip \citep{SaelensEtAl2003-EnvironmentalCorrelatesWalking}.
	\citet{HeinenEtal2010-CommutingByBicycle} study the determinants for bicycle commuting and find that the built environment affects traveler choice, among which distance is probably the most important factor.
	A study for British cities and towns by \citet{CerveroEtAl2019-NetworkDesignBuilt} suggests that safe connections that are as close to the shortest path as possible are most likely to encourage bicycle commuting.
	\citet{OspinaEtal2020-ColumbiaTravelDistance} arrive at a similar conclusion for Medellin city in Colombia: cyclists are willing to take a detour to ride on dedicated lanes, but only up to an extent.
	\citet{WangEtAl2021-DiffusionPublicBicycle} found that similar factors affect the choice to use shared bicycles.
	
	This paper presents optimization models to find the best way to improve an existing bicycle network for a given budget. Policy makers may use this method to guide their investments in cycling infrastructure, and to obtain the advantages that come with it. The optimization problem is formalized as the Bicycle Network Improvement Problem (BNIP) that selects how best to spend a given budget for road improvement in order to minimize the total penalty for a set of traveler origin-destination pairs (ODs).
	The penalties are calculated from the distance deviations from travelers' shortest paths.
	The optimization problem uses piecewise linear penalty functions, which are flexible enough to model different human preferences and numerous other factors.
	
	This paper models the BNIP as a Mixed-Integer Linear Program (MIP), which is solved through Benders decomposition.
	The optimization method is then used to conduct an in-depth case study for the Midtown area in Atlanta, GA, based on real transportation data.
	As the 10th most congested city in the US, and with a bicycle infrastructure scored in the red category \citep{Reed2019-GlobalTrafficScorecard}, Atlanta makes for an interesting test case. 
	Compared to the city of Delft in the Netherlands, cyclists in Atlanta are over two times more likely (78\% versus 32\%) to report poor road infrastructure as a cause of stress \citep{GadsbyEtAl2021-InternationalComparisonSelf}.
	Furthermore, a survey by the Atlanta Department of City Planning mentions that 70\% of people in the city currently feel uncomfortable to ride a bike \citep{Bottoms2018-CityAtlanta2018}.
	
	The paper contains four main contributions:
	\begin{enumerate}
		\item From a methodology standpoint, the paper formalizes the Bicycle Network Improvement Problem (BNIP) and shows that Benders decomposition is able to find optimal improvement plans for realistic instances, while the problem is computationally intractable for state-of-the-art black-box solvers.
		\item From a case study standpoint, the paper shows that even small investments in infrastructure may allow many additional commuters to travel safely by bike.
		\item At the intersection of methodology and case study, the paper demonstrates the value of optimization, which produces improvement plans that are significantly better than those obtained by heuristics. Moreover, and this is important for cities, the paper compares the benefits of optimal long-term plans with those obtained by upgrading the infrastructure incrementally. The paper shows that, on the case study, successive improvements using the BNIP lead to a network that is very close to optimal in the long term, which simplifies decision making.
		\item From a computational perspective, this paper compares a wide range of different penalty functions, and reports consistent results across all of them.
	\end{enumerate}
	
	\noindent
	The rest of the paper is organized as follows. Section \ref{sec:lit review}
	reviews prior work.
	Section \ref{sec:BNIP} formally introduces the BNIP, and Section \ref{sec:solution methods} describes a Benders decomposition algorithm to solve it. 
	Section \ref{sec:current atl} discusses the current conditions in Midtown Atlanta, which also motivates this research.
	The Midtown Atlanta case study is presented in Section \ref{sec:case study atl}.
	Section~\ref{sec:comparing BNIP formulations} explores the use of different penalty functions, and Section~\ref{sec:discussion} discusses how the methods in this paper can be adapted for future work. Finally, Section~\ref{sec:conclusion} presents the conclusions and possible directions for future research.

	\section{Review of Prior Work}
	\label{sec:lit review}
	There are several studies that consider bicycle infrastructure improvement planning.
	\citet{DuthieUnnikrishnan2014-OptimizationFrameworkBicycle} present a network design formulation to connect all OD (origin-destination) pairs with a lower bound on the bicycle level of service and an upper bound on the maximum travel length expressed as a function of the corresponding shortest path.
	Their objective and the BNIP objective are similar in that they limit the worst service for travelers with respect to the travel distances.
	There is, however, a fundamental difference between their work and the BNIP: the former mandates that all ODs admit feasible bicycle travels regardless of the improvement budget, where the BNIP has a limited budget to serve as many OD pairs as possible.
	The benefit of having a finite budget is that it abides by realistic scenarios, e.g., urban planners developing new bicycle infrastructure.
	Indeed, budgets for infrastructure improvements are almost always limited and their effective use is a key aspect for decision makers.
	
	\citet{MauttoneEtAl2017-BicycleNetworkDesign} introduce another MIP model to minimize the overall travel cost of riders, where the cost primarily consists of travel distances.
	This formulation includes a budget constraint, but still requires all OD pairs to be served.
	This is made possible by allowing for bicycle trips that are not 100\% safe, and penalizing the usage of unsafe roads.
	As safety is the primary concern for many potential cyclists, as argued in the introduction, the BNIP does not sacrifice the requirement for completely safe bicycle routes; rather it imposes limits on maximum travel distances to model realistic trips and provides an outside option for those riders who do not have a realistic safe route.
	It is also important to mention that \citet{MauttoneEtAl2017-BicycleNetworkDesign} only provide 
	sub-optimal solutions in reasonable time for their real-life case studies, and they use a heuristic to report results for large cases with more accuracy.
	The Benders decomposition algorithm proposed in this paper, however, solves the large Atlanta instances to optimality.
	
	\citet{LiuEtAl2020-UrbanBikeLane} present a MIP model to plan bicycle networks using objectives for coverage and continuity of travels. They assume that the MIP model receives, as input, bicycle paths. 
	Their adjacency-continuity utility function, which incorporates both safety and trip length, selects one of the pre-calculated paths to route each traveler while maximizing the utility of the network. 
	Their work is similar to the BNIP as it improves both safety and proximity of the trips. 
	The BNIP, however, has full flexibility in routing cyclists; this simplifies modeling for decision makers and may produce solutions of higher quality since the optimization can choose the best routes for riders and is not constrained by pre-selected paths.
	
	In addition to the previous works, a number of studies incorporate more diverse characteristics in the problem modeling.
	For example, \citet{LinYu2013-BikewayNetworkDesign}, \citet{LinLiao2016-SustainabilitySiBikeway}, and \citet{ZhuZhu2020-MultiObjectiveBike} use multi-objective optimization to include various objectives such as road connections, accessibility, and service level. 
	These formulations can model more customized bicycle experiences, but become more computationally intensive.
	The BNIP is a single-objective optimization problem with an objective that is flexible enough to model realistic applications.
	Most importantly, the proposed Benders decomposition algorithm is capable of performing studies in much larger instances and areas than those multi-objective programs. 
	It is an interesting avenue for future research to study if the techniques in this paper can be generalized to multiple objectives. 
	
	A number of other studies rely on heuristic methods instead of mathematical optimization techniques.
	\citet{BaoEtAl2017-PlanningBikeLanes} use large-scale bicycle trajectory data to define a flexible objective that combines the population covered by the network and the distances of their trajectories. It is solved with greedy-based heuristics that include steps to initiate road segments, expand the network, and terminate the improvement when the budget limit is met.
	Also, \citet{HsuLin2011-ModelPlanningBicycle} utilize shortest paths of ODs like those of \citet{DuthieUnnikrishnan2014-OptimizationFrameworkBicycle} to evaluate the quality of the network.
	They use a variety of algorithms, some of which are greedy, to compare the shortest paths to bicycle routes.
	\citet{OrozcoEtal2020-DataDrivenOptimalNetwork} introduce two greedy algorithms to connect bicycle network components.
	They also compare the shortest paths to bicycle routes, but as opposed to the BNIP, they do not optimize any objective.

	\section{The Bicycle Network Improvement Problem} 
	\label{sec:BNIP}
	This section introduces the Bicycle Network Improvement Problem (BNIP) to find the best improvement of an existing bicycle network within a given budget.
	Let the current road network be represented by a directed graph $G=(V,W)$ with nodes $V$ and arcs $W$. 
	The arcs are referred to as \emph{ways}. 
	Ways are partitioned into two distinct sets $W = W^{safe} \cup W'$, with $W^{safe}$ the set of ways that are safe for cycling, and $W'$ the set of unsafe ways. Every way $(i,j) \in W$ has a length, given by the parameter $d_{ij} \ge 0$.
	The total length of ways that can be improved is limited by the budget $B$.
	The set of sample trips that travel through the network is given by $T$.
	Each OD $k \in T$ consists of an origin $o_k \in V$ and a destination $d_k \in V$. Additionally, let $s_k \ge 0$ be the length of the (possibly unsafe) shortest path between $o_k$ and $d_k$, and let $p_k \ge 1$ be the number of travelers completing this travel.

	\subsection{Modeling Bicycle Travel}
	\label{sec:characteristics}
	
	Safety is critical to increase cyclist participation. Accordingly, two
	characteristics are taken into account when modeling bicycle trips.
	First, potential cyclists would like completely safe routes from
	origin to destination: if the safety requirement is met for a certain
	OD, then the route is labeled as \emph{safe}. If the network cannot
	provide a safe route, the BNIP assumes that the potential rider will
	select another mode of transportation, which is referred to as the
	\emph{outside option}.  Second, the travel should not take much longer
	than the alternative transportation mode, e.g., driving by car. Hence,
	if there is no safe path of length smaller than $L_k$ (a parameter for
	rider $k$), the BNIP assumes that rider $k$ will not travel by
	bicycle.
	
	To model the appeal of short bicycle trips for rider $k$, the BNIP
	uses a penalty function $f_k$ that is non-decreasing and satisfies
	$f_k(0) = 0.$ For an OD pair $k$ and a path of length $l_k$, the
	penalty is given by $f_k(u_k)$ where $u_k = l_k - s_k$ denotes the
	deviation from the (potentially unsafe) shortest path.  If rider $k$
	cannot be provided a short trip, the outside option is used, and the
	penalty is the objective function is defined as $L_k -
	s_k$ (An optimization model that purely maximizes the number
	of cyclists is also presented in the paper).  As such, this trip is
	assigned the same penalty as a path of length $L_k$, which is the
	tipping point at which the rider starts preferring the outside option.
	
	\subsection{Mathematical Formulation}
	\label{sec:formulation for BNIP-P}
	The mathematical model for the BNIP is depicted as follows: 
	\begin{spacing}{1}
		\begin{mini!}
			%
			{}
			%
			{
				\sum_{k \in T} p_k \ f_k(u_k),
				\label{eq:new:objective}
			}
			%
			{\label{formulation:new:BNIP}}
			%
			{}
			%
			%
			\addConstraint
			{
				\sum_{(i,j)\in W'}
			}
			{d_{ij} y_{ij} \leq B,\hspace{4.4cm}}
			{\label{eq:new:budget}}
			\addConstraint
			{
				\hspace{-1.1cm}
				\mathrlap{
					\sum_{(i,j)\in W} x^k_{ij} - \sum_{(j,i) \in W} x^k_{ji}
					= \begin{cases}
						1 - z_k & \textrm{ if } i = o_k\\
						z_k - 1 & \textrm{ if } i = d_k \\
						0 & \textrm{ otherwise}
					\end{cases}
				}
			}
			{}
			{\,\,\,\,\,\,\,\,\,\,\,
				\forall k\in T, i\in V,
				\label{eq:new:flow_balance}
			}
			\addConstraint
			{x^k_{ij}}
			{\leq y_{ij}}
			{\forall k \in T, (i,j)\in W', \label{eq:new:upgraded_way_avail}}
			\addConstraint
			{u_k}
			{\ge \sum_{(i,j) \in W} d_{ij} x^k_{ij} +  L_k z_k - s_k}
			{\forall k \in T, \label{eq:new:deviation}}
			\addConstraint
			{y_{ij}}
			{\in \mathbb{B}}
			{\forall (i,j)\in W', \label{eq:new:s_binary}}
			\addConstraint
			{x^k_{ij}}
			{\in \mathbb{B}}
			{\forall k \in T, (i,j)\in W, \label{eq:new:x_binary}}
			\addConstraint
			{z_{k} }
			{\in \mathbb{B}}
			{\forall k \in T. \label{eq:new:z_binary}}
		\end{mini!}%
		\doublespacing
	\end{spacing}
	For every way $(i,j)\in W'$, variable
	$y_{ij} \in \mathbb{B}$ indicates whether $(i,j)$ is upgraded to safe
	conditions (value one), or remains unchanged (value zero). The
	shortest safe path for every OD is determined by variables $x$ and
	$z$: variable $x^k_{ij} \in \mathbb{B}$ represents whether trip $k\in
	T$ uses way $(i,j)$. Variable $z_k \in \mathbb{B}$ indicates whether
	trip $k \in T$ uses the outside option.  As explained in the previous
	section, the variable $u_k$ represents the argument of the penalty
	function for every trip $k \in T$.
	
	Objective~\eqref{eq:new:objective} minimizes the total penalty of the
	riders $\sum_{k \in T} p_k f_k(u_k)$ on the network. The penalties for
	cyclists are computed in terms of the deviation from the shortest-path
	distance. Riders who are not cycling or have safe paths that are too
	long incur the penalty associated with a path of length $l_k=L_k$, as
	will become clear shortly. Constraint~\eqref{eq:new:budget} limits the
	budget for improving the network.
	Constraints~\eqref{eq:new:flow_balance} impose the path conservation
	conditions: each OD has either a unit flow (if $z_k=0$), in which case
	the $x$-variable describes the path, or uses the outside option
	($z_k=1$). Constraints~\eqref{eq:new:upgraded_way_avail} make sure
	that unsafe ways can only be used if
	upgraded. Constraints~\eqref{eq:new:deviation} compute the deviation.
	Note that, if the shortest path for trip $k$ exceeds length $L_k$, it
	is optimal to set $z_k$ to 1, i.e., trip $k$ uses the outside option.
	Infeasible trips are also assigned the penalty associated with a path
	of length
	$L_k$. Constraints~\eqref{eq:new:s_binary}-\eqref{eq:new:z_binary}
	capture the integrality conditions.

	\section{Solution Methods for the BNIP} 
	\label{sec:solution methods}
	
	Solving the BNIP directly with a MIP solver, such as CPLEX or Gurobi,
	is computationally intractable for the scale of the case study
	considered in this paper.  Observe however that, for a given design
	(i.e., when the $y$-variables are fixed), the formulation reduces
	to a set of independent minimum-cost flow problems, one for each
	OD. By total unimodularity, this implies that the integrality
	conditions \eqref{eq:new:x_binary}-\eqref{eq:new:z_binary} can be
	relaxed. This makes the problem ideally suited for Benders
	decomposition \citep{Benders1962-PartitioningProceduresSolving}.
	
	\subsection{Benders Decomposition} 
	\label{sec:benders}
	
	The Benders decomposition for the BNIP has a master problem to determine how to upgrade the network, and subproblems to return the safe path or the outside option for each rider for a given upgraded network. 
	The master problem generates a network design, and the subproblems find the paths in the proposed network for each trip.
	The optimal solutions of the subproblems are then used to derive the Benders cuts that are added to the master problem.
	These two steps are iterated until no more violated Benders cuts are generated by the subproblems, at which point the network is optimal.
	
	\paragraph{The Benders Master Problem}
	The master problem is presented as follows: 
	\begin{mini!}
		%
		{}
		%
		{
			\sum_{k \in T} p_k \ f_k(u_k),
			\label{eq:new:master_obj}
		}
		%
		{\label{formulation:new:master problem}}
		%
		{}
		%
		%
		\addConstraint
		{
			\eqref{eq:new:budget}, \eqref{eq:new:s_binary}, \notag
		}
		{}
		{}
		\addConstraint
		{u_{k} \ge \Phi_k(y) \quad \forall \, k \in T. \label{eq:new:master_cut}}
		{}
		{}
	\end{mini!}%
	Its Objective~\eqref{eq:new:master_obj} is the same objective as in Formulation~\eqref{formulation:new:BNIP}.
	Constraints~\eqref{eq:new:budget} and \eqref{eq:new:s_binary} ensure that the network improvement plan is within budget and valid. 
	When solving the master problem, Constraints~\eqref{eq:new:master_cut}, where $\Phi_k(y)$ is the minimum objective value for trip $k$ given a design $y$, are replaced by the Benders cuts generated from the solutions to the subproblems. 
	
	\paragraph{Benders Subproblem}
	\label{sec:sub problem}
	The subproblem for a trip $k$ generates Benders cuts for each network
	produced by the master problem. The $z$-variables (the outside option) ensure \emph{complete recourse}: the subproblem is feasible for any network, because ODs can always use the outside option. This implies that the optimality cuts \eqref{eq:new:master_cut} are sufficient and no feasibility cuts are needed.
	
	For a given network, the subproblem decomposes into many independent subproblems, and the subproblem of the BNIP for each $k\in T$ is defined as follows:
	\begin{alignat}{4}
		\Phi_k(y) = \min \,&\,\,  \sum_{(i,j) \in W} d_{ij} x^k_{ij} + L_k z_k  - s_k, & & & & \\
		\text{s.t.} \,\,
		& \eqref{eq:new:flow_balance}, \eqref{eq:new:upgraded_way_avail},
		\eqref{eq:new:x_binary}, \eqref{eq:new:z_binary}. \,\,\, \notag &&
	\end{alignat}
	
	Since the subproblem is a standard minimum-cost flow problem, it is totally unimodular and can be solved by linear programming. 
	This implies that the Benders subproblem can be solved by optimizing many small and independent linear programs, which is the prime reason why the Benders decomposition provides a significant computational benefit.
	
	With dual variables $\lambda$ and $\mu$ associated with Constraints~\eqref{eq:new:flow_balance} and \eqref{eq:new:upgraded_way_avail}, respectively, the dual subproblem is defined as follows:
	\begin{maxi!}
		%
		{}
		%
		{
			\lambda_{o_k}^k - \lambda_{d_k}^k - \sum_{(i,j)\in W'} \mu_{ij}^k y_{ij} - s_k,
			\label{eq:sub_obj}
		}
		%
		{\label{formulation:sub problem dual}}
		%
		{\Phi_k(y) = }
		%
		%
		\addConstraint
		{\lambda_i^k - \lambda_j^k \le d_{ij}}
		{}
		{\forall \, (i,j) \in W^{safe}, \label{eq:sub_1}}
		\addConstraint
		{\lambda_i^k - \lambda_j^k - \mu_{ij}^k \le d_{ij}\qquad}
		{}
		{\forall \, (i,j) \in W', \label{eq:sub_2}}
		\addConstraint
		{\lambda_{o_k}^k - \lambda_{d_k}^k \le L_k,}
		{}
		{\label{eq:sub_3}}
		\addConstraint
		{\lambda_i^k \in \mathbb{R}}
		{}
		{\forall \, i \in V,\label{eq:sub_4}}
		\addConstraint
		{\mu_{ij}^k \ge 0}
		{}
		{\forall \, (i,j) \in W'. \label{eq:sub_5}}
	\end{maxi!}%
	
	The Benders cuts are obtained as
	\begin{equation}\label{eq:bender cuts}
		u_k \ge \bar{\lambda}_{o_k}^k - \bar{\lambda}_{d_k}^k - \sum_{(i,j)\in W'} \bar{\mu}_{ij}^k y_{ij} - s_k,
	\end{equation}
	where $\bar{\lambda}$ and $\bar{\mu}$ are optimal dual values for the corresponding trips.
	
	\paragraph{Pareto-Optimal Cuts}
	
	It is well-known that network flow problems often suffer from dual
	degeneracy.  \citet{MagnantiWong1981-AcceleratingBendersDecomposition}
	addressed this issue by generating \emph{Pareto-optimal cuts} that are
	not dominated by any other Benders cut.  This requires solving a
	Pareto subproblem that uses the result from the standard
	subproblem. Pareto-optimal cuts need a \emph{core point}, i.e., a
	point in the relative interior of the feasible region of the master
	variables.  For the BNIP, the following point $y'$ is selected
	as the core point.
	\begin{equation}
		y'_{ij} = \frac{1}{2} \min\left\{\frac{B}{\lvert W' \rvert d_{ij}}, 1\right\} \, \forall\, (i,j) \in W'
	\end{equation}
	The point $y'$ is in the relative interior, as $y'_{ij} \in (0,1)$ and $\sum_{(i,j) \in W'}  d_{ij} y'_{ij} \le \sum_{(i,j) \in W'}  d_{ij} \frac{B}{2\lvert W' \rvert d_{ij}} = \frac{B}{2} < B$, which strictly satisfies budget constraint~\eqref{eq:new:budget}.
	
	Using the core point $y'$, the Pareto subproblem is defined as follows:
	\begin{maxi!}
		%
		{}
		%
		{
			\lambda_{o_k}^k - \lambda_{d_k}^k - \sum_{(i,j)\in W'} \mu_{ij}^k y'_{ij} - s_k,
			\label{eq:pareto objective}
		}
		%
		{\label{formulation:pareto}}
		%
		{}
		%
		%
		\addConstraint
		{\lambda_{o_k}^k - \lambda_{d_k}^k - \sum_{(i,j) \in W'} \mu_{ij}^k y_{ij} - s_k = \Phi_k(y),}
		{}
		{}
		\addConstraint
		{\eqref{eq:sub_1}-\eqref{eq:sub_5}. \notag}
		{}
		{}
	\end{maxi!}%
	To use Pareto-optimal cuts, each Benders iteration is changed as follows.
	For every trip, the value of $\Phi_k(y)$ is calculated by solving the subproblem.
	Next, the Pareto subproblem is solved to produce new optimal dual values $\bar{\lambda}$ and $\bar{\mu}$, and those variables are used to generate cuts as in Equation~\eqref{eq:bender cuts}.
	
	\paragraph{Two-Phase Benders}
	\citet{McDanielDevine1977-ModifiedBendersPartitioning} observe that it
	is not necessary to solve the master problem to optimality at every
	iteration to obtain valid Benders cuts.  They propose to apply Benders
	decomposition in two phases.  In phase one, Benders decomposition is
	applied to the relaxed master problem.  For the BNIP, this amounts to
	relaxing the integrality conditions~\eqref{eq:new:s_binary} in the
	master problem, and solving the subproblems for fractional values of
	$y$.  In phase two, the integrality conditions are reinstated,
	and the Benders decomposition algorithm continues with the original
	master problem.  The Benders cuts that are added in phase one are
	maintained, which ensures a better starting lower bound, which often
	improves the overall performance of the algorithm.  Moving from phase
	one to phase two is possible at any point, and this paper switches
	over when the relaxed problem is solved, or when a time limit is
	reached.

	\subsection{Greedy Heuristic}
	\label{sec:heurstic}
	This section also introduces a greedy heuristic that only relies on the ability to solve shortest paths to demonstrate the value of optimization.
	The heuristic greedily computes the next way to upgrade in the network. 
	For every OD, it computes a shortest path that minimizes the total distance on unsafe roads, i.e., the costs for traveling safe or unsafe way $(i,j)\in W$ is 0 or $d_{ij}$ respectively.
	Only the shortest paths with a distance of at most $L_k$ are considered for all $k \in T$. 
	The relative importance of an unsafe way is determined by counting the total number of riders whose shortest paths include this way.
	The greedy heuristic selects the most ``important" way to upgrade and repeats the process until the budget is exhausted.

	\section{Cycling in Midtown Atlanta}
	\label{sec:current atl}
	This paper is motivated by improving bicycle travel conditions in the Midtown area in Atlanta, GA.
	Midtown is a neighborhood of Atlanta that consists of a commercial core and a residential neighborhood.
	The Midtown core is characterized by high-rise buildings and functions as a major employment center, including offices of large companies such as NCR, Google, Equifax, and Honeywell.
	The commercial core has around 28k residents, and 70k people per day travel to the area for work \citep{MidtownAlliance2019-2019MidtownCommunity}.
	The residential neighborhood is to the east of the business district, and mainly consists of single-family residences.
	Every workday, a large number of commuters drive to Midtown, and cause a significant amount of traffic congestion.
	To better this situation, the case study aims \emph{at improving the bicycle infrastructure to provide these commuters safe and short cycling trips} as an alternative to commuting by car.
	
	\subsection{Travel Data}
	To gain insight into current commuting behaviors, this paper uses travel data provided by the Atlanta Regional Commission (ARC).
	The ARC used an activity-based model, calibrated with survey data collected in 2007-2011, to simulate trips at the individual level \citep{ARC2017-ActivityBasedModeling}.
	The data reveals a low share of cyclists among commuters to Midtown (about 0.7\%), but also a significant potential for improvement: many of the trips to Midtown are short, and over 70\% of commutes are completed by people driving alone, who could potentially switch to cycling if the infrastructure is improved.
	
\begin{figure}[!t]
		\centering
		\includegraphics[width=0.7\textwidth]{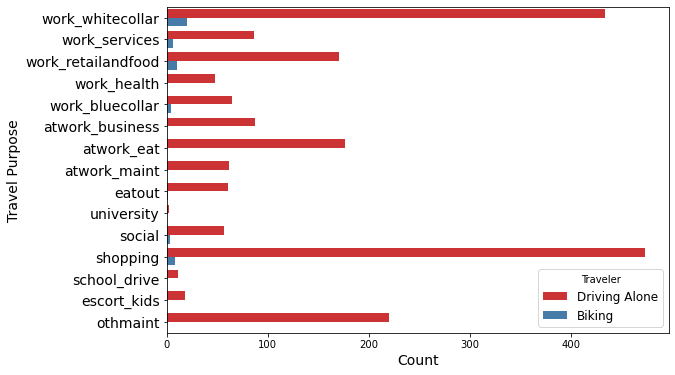}
		\caption{Travelers from Virginia-Highland to Midtown.}
		\label{fig:abm_trend}
	\end{figure}

	The case study focuses on one particularly interesting group: the group of white-collar workers coming in from \emph{Virginia-Highland}, an affluent neighborhood immediately to the east of Midtown.
	Figure~\ref{fig:abm_trend} shows the daily number of trips originating from Virginia-Highland for both cyclists and solo drivers, categorized by trip purpose.
	Many trips are taken by white-collar commuters but, despite the close proximity of Virginia-Highland to Midtown, the number of cyclists is less than 5\% that of the number of solo drivers.
	
\begin{figure}[!t]
		\begin{minipage}[t][][b]{0.5\textwidth}
			\centering
			\includegraphics[width=0.8\linewidth]{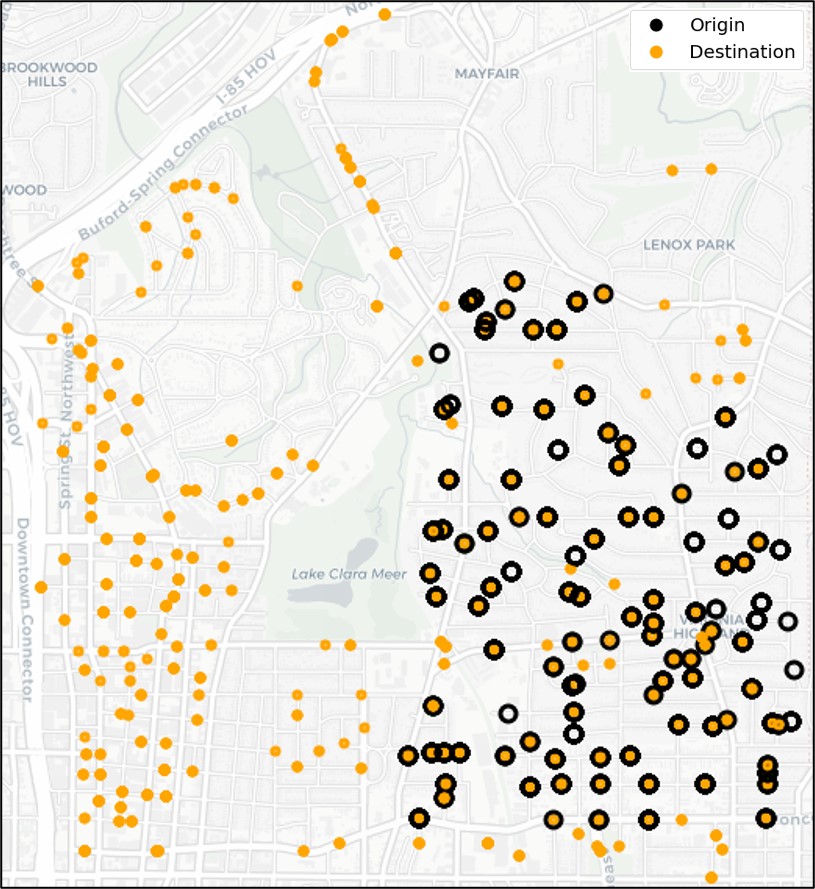}
			\caption{Origins and Destinations of 1039 sample OD pairs.}
			\label{fig:sample_od}
		\end{minipage}
		\begin{minipage}[t][][b]{0.5\textwidth}
			\centering
			\includegraphics[width=0.8\linewidth]{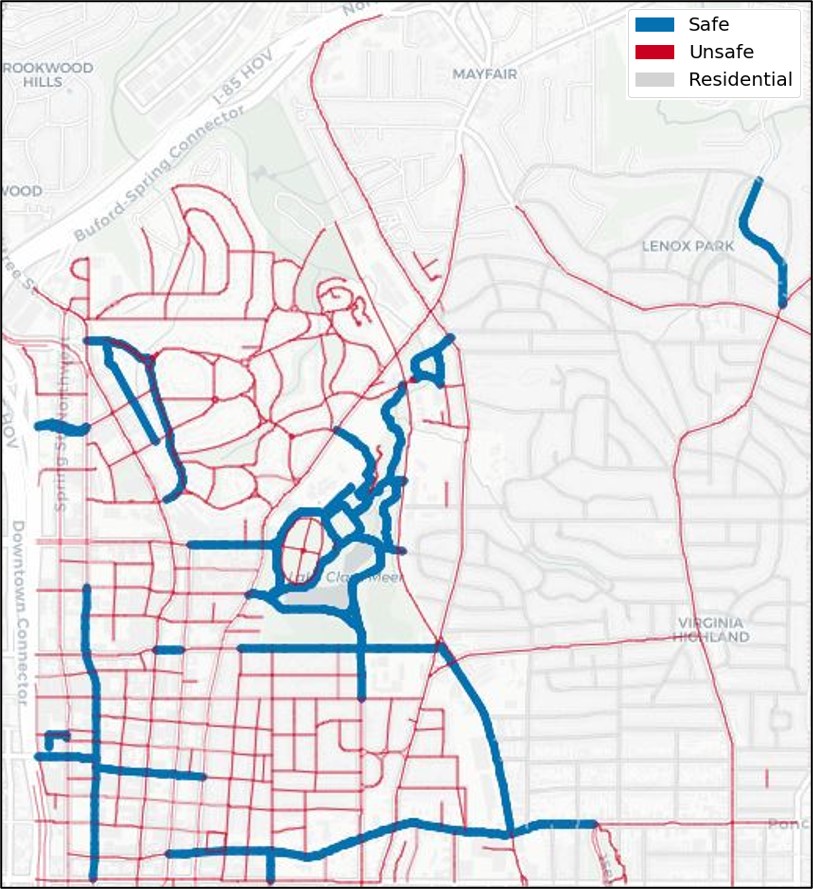}
			\caption{Current Bicycle Network.}
			\label{fig:current}
		\end{minipage}
	\end{figure}

	To study how improving the bicycle network affects travelers, a sample of travels is generated to represent the white-collar commuters from Virginia-Highland.
	The eight Traffic Analysis Zones (TAZs) that cover Virginia-Highland are selected as the departure zones, and 72 TAZs that cover Midtown and Virginia-Highland are chosen as the destination zones.
	Destinations in Virginia-Highland are included, as these trips may benefit from the same infrastructure improvements.
	The ARC provides travel data between the TAZs, and to obtain a more realistic sample, the origins and destinations are randomly assigned to the centers of the smaller census blocks, weighted by the population counts.
	Five samples are generated for every OD pair, and samples that do not connect to the existent road network are filtered out.
	The result is a set of 1039 representative trips, covering 110 origins and 256 destinations, presented by Figure \ref{fig:sample_od}.
	The origins on the map are in Virginia-Highland in the East, and most destinations are in the business area in the West.
	The center of the map shows Piedmont Park, with the Midtown residential neighborhood to its South.
	
	\subsection{Current Bicycle Network}
	\label{sec:current atl situation}
	The current road network in the case study area was retrieved from OpenStreetMap \citep{OpenStreetMap2020}.
	The network consists of 5815 nodes and 11,329 directed ways.
	The ways make up 667 roads, and have a total length of 339 km (212 miles).
	The roads are classified into three types: roads with a dedicated bike lane, residential roads, and unsafe roads, where dedicated and residential roads are assumed to be safe for cyclists.
	
	Figure \ref{fig:current} shows the current bicycle network in the case
	study region, where blue, red, and gray roads describe dedicated,
	unsafe, and residential roads, respectively.  Unsafe roads do not have
	the proper infrastructure for cyclists, and are the target for
	conversion to dedicated roads.  There are 450 bicycle-unsafe roads,
	with a total length of 180 km (113 miles).  The conditions on the unsafe roads
	are assumed to be similar throughout the area, and the cost to realize
	dedicated bike lanes is assumed to be the same per unit length
	everywhere. Among the sampled ODs, only 170 trips (16\%) have 
	bicycle-safe routes that are completely safe. 
	Moreover, many of them require a significant detour
	to complete those safe paths.  For example, only 89 ODs (9\%) have
	access to a bicycle path with a detour of less than 10\% of the
	shortest path.

	\section{Experimental Results}
	\label{sec:case study atl}
	
	This section presents experimental results on the case study. It first shows that the Benders decomposition algorithm allows for solving the BNIP to optimality for realistic instances that cover a whole neighborhood. It also analyzes the effectiveness of the optimal improvement plans to provide safe and short cycling trips to commuters
	and demonstrates the benefits of optimization by comparing the optimal improvement plans to those obtained with the greedy heuristic.
	
	\subsection{Experimental Settings}
	\label{sec:experimental settings}
	The experimental results use the case study region, the existing
	network, and the 1039 sample trips from Section~\ref{sec:current atl}.
	For each sample, the number of passengers $p_k=1$, as the case study
	targets solo drivers that commute to their work places. The budget $B$
	ranges from 6.4 km (4 miles) up to 44.8 km (28 miles) of improvements, in 6.4 km (4 mile)
	increments.  For the case study, it is useful to consider network
	improvements on the road level, rather than on the way level, and to
	improve both directions at the same time.  These solutions are more
	practical to implement, and also contribute to safe return trips.  In
	total, 450 roads are identified.  To improve all ways of a road at the
	same time, additional constraints are added to the BNIP: if $(i,j)\in
	W'$ and $(p,q)\in W'$ are part of the same road, then $y_{ij} =
	y_{pq}$.  These constraints are added to the Benders master problem,
	without changing the main steps of the algorithm. For the greedy
	heuristic, the relative importance of a road is calculated by averaging the
	usage counts for the individual ways. 
	
	The core experimental results use a linear penalty function with $f_k(u_k)=u_k$, i.e., the objective maximizes the number of cyclists and minimizes the overall average distance over the shortest paths, and the resulting BNIP is denoted as BNIP-L.
	Section~\ref{sec:comparing BNIP formulations} provides experiments with alternative penalty functions, and Section~\ref{sec:discussion} discusses how additional data can be included in the model to complement travel distance.
	The threshold $L_k$ is computed in terms of a deviation factor $R\ge 1$, i.e., $L_k=s_k R$.
	The full Benders algorithm is implemented in Python, and Gurobi 9.0.2 is used to solve the master problem and the subproblems.
	All computations are performed with an Intel Core
	i7-8565U CPU and 16GB of RAM.
	
	\subsection{Efficiency of the Benders Algorithm}
	\label{sec:efficiency}
	
	This section compares three Benders decomposition algorithms: the
	traditional Benders decomposition algorithm (TB), the algorithm that
	uses \citeauthor{MagnantiWong1981-AcceleratingBendersDecomposition}
	Pareto-optimal optimal cuts (MW), and the algorithm that uses both
	Pareto-optimal cuts and the \citeauthor{McDanielDevine1977-ModifiedBendersPartitioning} two-phase
	strategy (MW-McD).  In the MW-McD experiments, phase one is limited to
	20 minutes.
	\begin{figure}[!t]
		\begin{adjustbox}{center}
			\includegraphics[width=\linewidth]{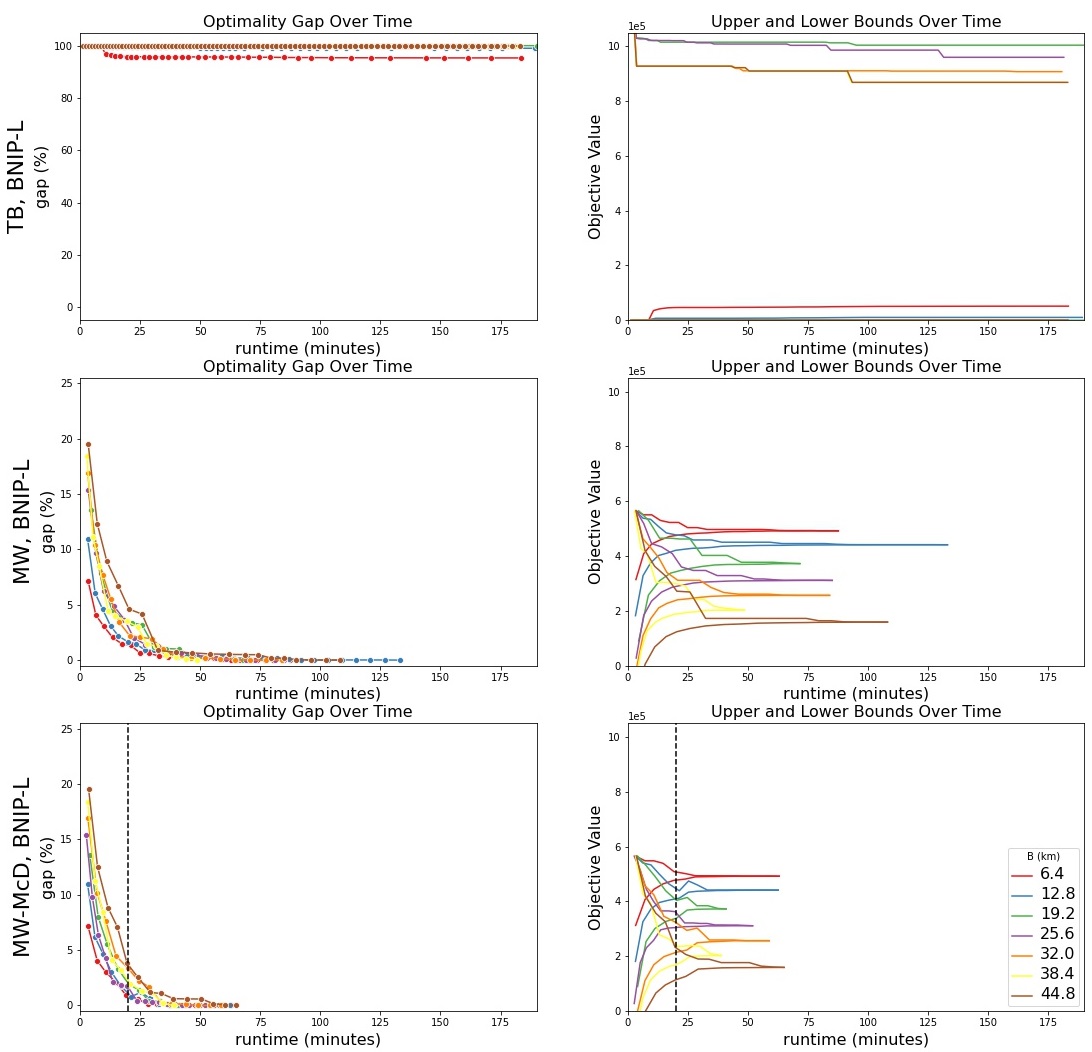}
		\end{adjustbox}
		\caption{The Performance of Benders Decomposition Algorithms ($R=1.2$).}
		\label{fig:gap}
	\end{figure}
	
	Figure~\ref{fig:gap} reports the experimental results for various budgets and $R=1.2$.
	The left charts show the optimality gap, i.e., the difference between the best feasible solution and the lower bound over time, and the right charts show how the upper and lower bounds approach each other until the optimal solution is found.
	The three rows correspond to the TB, MW, and MW-McD strategy, respectively.
	In the MW-McD case, the vertical dashed line indicates the switch from phase one to phase two.
	
	\begin{spacing}{1}
		\begin{table}[!t]
			\centering
			\begin{tabular}{c r r}
				\toprule
				& \multicolumn{2}{c}{Iterations Count} \\
				\cmidrule{2-3}
				B (km) & MW & MW-McD\\
				\midrule
				6.4 & 24 & 6+10=16\\
				12.8 & 28 & 6+10=16\\
				19.2 & 16 & 5+5=10\\
				25.6 & 20 & 4+13=17\\
				32.0 & 17 & 6+7=13\\        
				38.4 & 13 & 6+4=10\\
				44.8 & 19 & 6+8=14\\
				\bottomrule
			\end{tabular}
			\caption{Number of Iterations for MW and MW-McD (Fractional+Integral Cuts).} \label{tab:benders iteration}
		\end{table}
		\doublespacing
	\end{spacing}
	
\begin{figure}[!t]
		\begin{adjustbox}{center}
			\includegraphics[width=\linewidth]{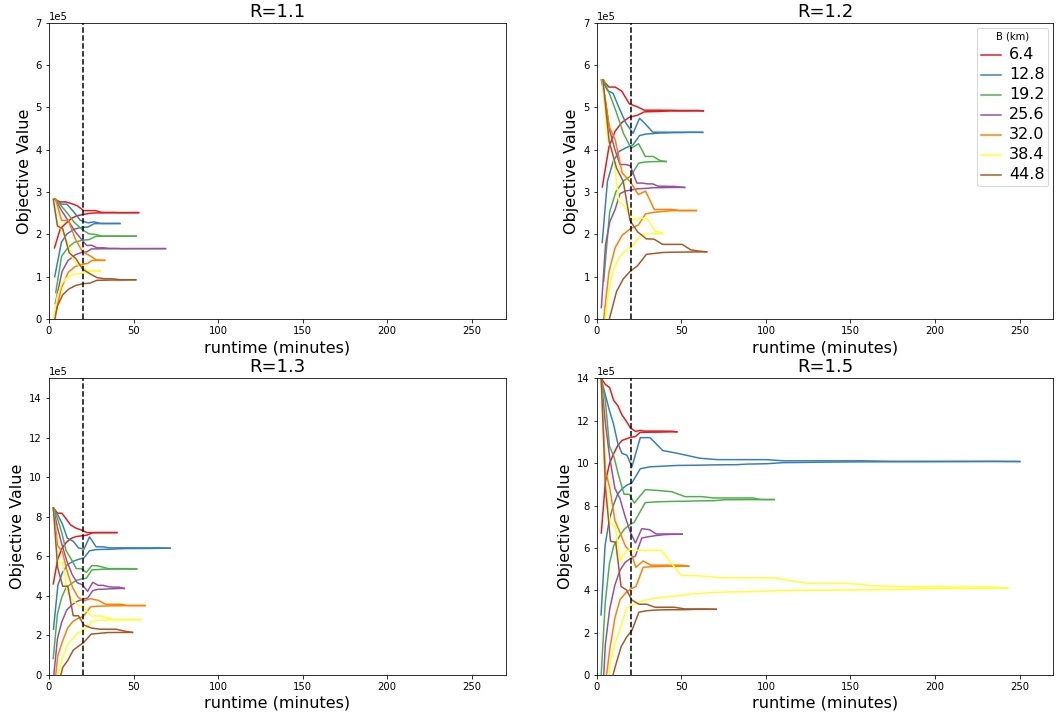}
		\end{adjustbox}
		\caption{The Performance of MW-McD Algorithm with Different Distance Thresholds.}
		\label{fig:gap different R}
	\end{figure}

	The first observation is that TB is significantly outperformed by MW
	and MW-McD.  When fractional cuts are added prior to integral cuts
	(MW-McD), the initial optimality gap is much smaller and the number of
	iterations is significantly reduced compared to MW.
	Table~\ref{tab:benders iteration} compares MW and MW-McD on the number
	of Benders iteration to reach optimality. MW solves each BNIP instance
	optimally in under 150 minutes, where MW-McD only takes
	around an hour per instance.  Figure~\ref{fig:gap different R}
	verifies that the good performance of MW-McD is consistent for
	different budgets $B$ and for different distance thresholds $R$.
	Overall, the solution time of at most five hours is short for creating
	an improvement plan for months or years into the future.
	This contrasts with the MIP model that cannot solve the instances in
	reasonable time.  Note that the difficulty comes from the size of the
	problem: the network consists of $\lvert W \rvert = 11,329$ arcs, and
	the BNIP introduces a flow variable for each arc and for each of the
	$1039$ ODs, resulting in a MIP with over 10 million variables.  The
	Benders decomposition exploits the problem structure and solves a
	significantly smaller problem for each OD at every step.
	
	\subsection{Impact of Bicycle Network Improvement Plans}
	\label{sec:road improvement result}
	It is interesting to study the type of improvements produced by
	optimal plans.
	Figures~\ref{fig:examplemaps}-\ref{fig:examplemaps_1.5} in
	Appendix~\ref{sec:appen_example maps} present the full series of
	improvement plans for the different settings of $B$ and $R$; this
	section discusses the most important observations.  The experiments
	are conducted with $R=1.1, 1.2, 1.3$, and $1.5$.

	\paragraph{Plan Characteristics}
	
	The bicycle network improvement plans show several notable trends.
	First, the problem searches for sub-regions that can be served with minimum improvements by maximizing the usage of pre-existing infrastructure. The left map of Figure~\ref{fig:4,12mi maps} shows the optimal plan when only 6.4 km (4 miles) of road can be improved.
	\begin{figure}[!t]
		\centering
		\includegraphics[width=0.8\textwidth]{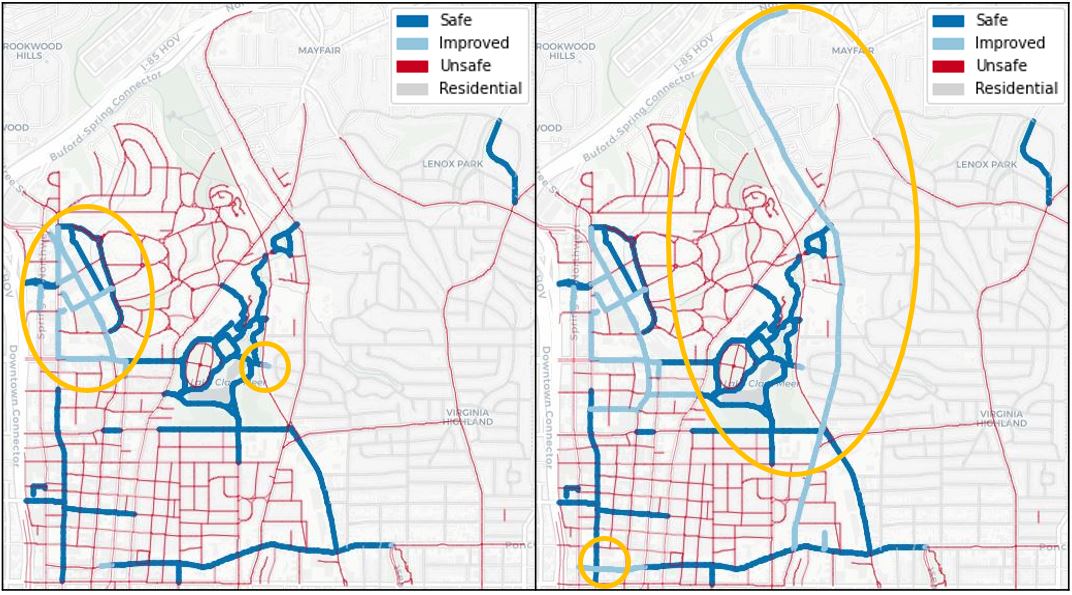}
		\caption{Optimal Improvement Plans ($R=1.1$): $B= 6.4$ km (4mi, left) and $19.2$ km (12mi, right).}
		\label{fig:4,12mi maps}
	\end{figure}
	Two crucial improvements (circled in yellow) were selected: a short segment that links the park and Virginia-Highland to the East, and roads Northwest of the park that connect to the business area. These connections are essential to allow the commuters from Virginia-Highland to commute to their workplaces.
	
	A second observation is that the optimal plans do not waste budget to improve multiple roads that serve a similar purpose. For instance, consider the short discontinuity of safe roads between the park and the roundabout located north of the park. The optimization algorithm does not remove this discontinuity, because the short segment connecting Virginia-Highland to the park serves the same purpose. Plans generated by the greedy heuristic (Figure~\ref{fig:examplemaps_HEU} in  Appendix~\ref{sec:appen_example maps}) do not recognize this. This is the value of optimization that provides globally optimal plans. This will have the consequence that some trips are better served on the heuristic network, but overall the optimization will produce significantly better plans. 
	
	A third observation concerns increased budgets: here the BNIP
	prioritizes the improvement of the backbone of the network, rather
	than providing sporadic developments. The 19.2 km (12 mi) map in Figure
	\ref{fig:4,12mi maps} shows the investments on Monroe Drive (a
	hook-like vertical road) to provide safe North-South travel, and on
	missing segments of Ponce De Leon Avenue (a horizontal road located
	southwest), that may connect all Southern demands. 
	
\begin{figure}[!t]
		\begin{minipage}[p][][b]{0.5\textwidth}
			\centering
			\includegraphics[width=0.95\textwidth]{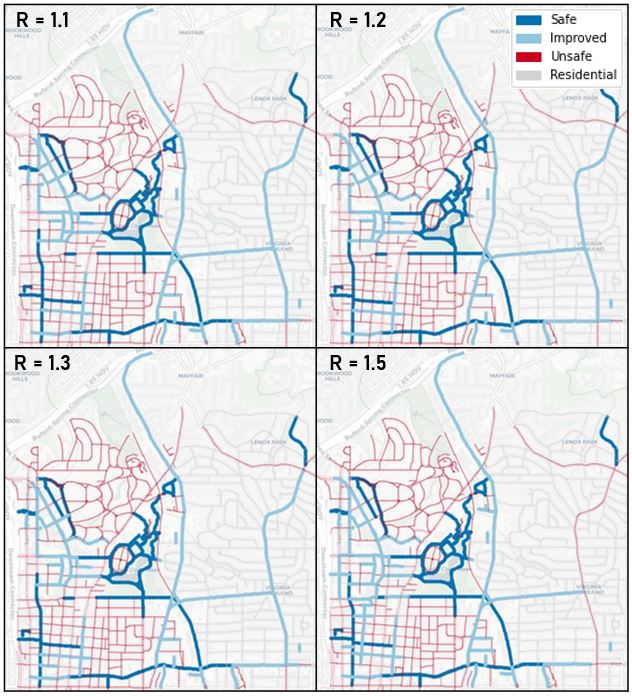}
			\caption{Optimal Improvement Plans, $B=38.4$ km (24 mi).}
			\label{fig:1.1 1.2 1.3 1.5 comparison}
		\end{minipage}
		\begin{minipage}[p][][b]{0.5\textwidth}
			\centering
			\includegraphics[width=0.95\textwidth]{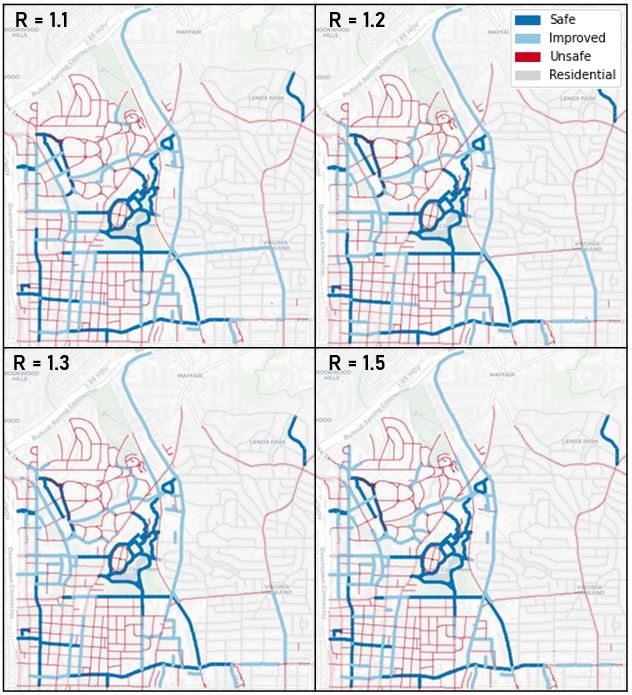}
			\caption{Heuristic Improvement Plans, $B=38.4$ km (24 mi).}
			\label{fig:1.1 1.2 1.3 1.5 comparison - heu}
		\end{minipage}
	\end{figure}

	The order in which the backbone is constructed depends on the deviation factor. 
	For instance, Figure \ref{fig:1.1 1.2 1.3 1.5 comparison},
	which includes 38.4 km (24 miles) improvement plans with four deviation factors, shows that
	North Highland Avenue (a vertical road located east) has not been
	improved for $R=1.5$, unlike three other plans of shorter distance
	thresholds. That is because the longer deviation allowance permits to
	serve origins located on the east with residential roads despite
	causing some detour, and improves more roads in the western business
	area to complete more last miles of the ODs. 
	
	Although the order of its construction may be different, the backbone of the network converges 
	as $B$ increases. To follow up on the above case, the 44.8 km (28 miles) plan for $R=1.5$
	(Figure~\ref{fig:examplemaps_1.5} in Appendix~\ref{sec:appen_example maps}) 
	improves North Highland Avenue and exhibits practically
	identical road improvements regardless of different travel length allowance.

\begin{figure}[!t]
		\centering
		\includegraphics[width=0.7\textwidth]{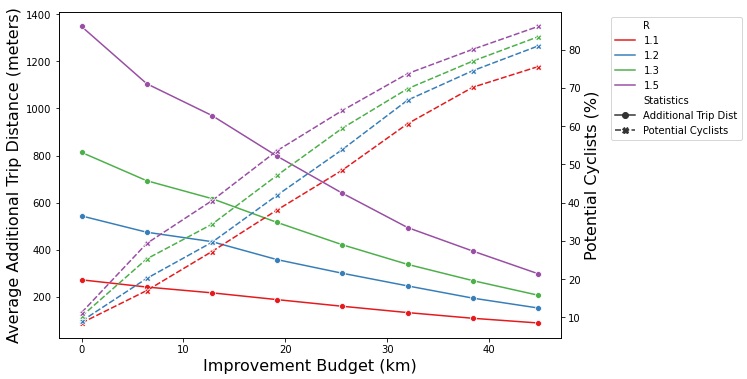}
		\caption{Effectiveness of the BNIP Improvements.}
		\label{fig:current_vs_opt}
	\end{figure}

	\paragraph{Effectiveness}
	Figure \ref{fig:current_vs_opt} shows the effectiveness of the optimal improvement plans for different parameters.
	For each $R$, the percentage of potential cyclists grows as the budget
	and the average trip distance over the shortest paths
	decreases. Potential cyclists are those riders with a safe bicycle
	trip whose length does not exceed the maximum distance (i.e, $L_k$ for
	trip $k$). The number of people that benefit from the improvements is
	similar for all values of $R$. Improving only 6.4 km of bicycle lanes
	already doubles the amount of potential cyclists at the
	minimum. Moreover, the number of potential cyclists increases almost
	linearly with the budget, suggesting that, on the case study, further investments in bicycle infrastructure deliver similar value and keep increasing the number of potential cyclists with a safe and short route.
	
\begin{figure}[!t]
		\centering
		\includegraphics[width=0.85\textwidth]{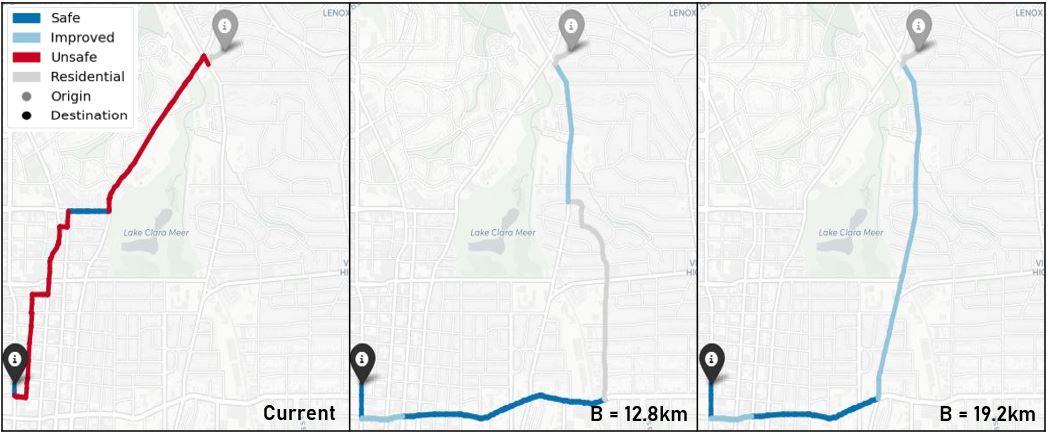}
		\caption{Example Routes on Optimal Improvement Plans, R=1.5.}
		\label{fig:orig_8_12mi_k973}
	\end{figure}

	Figure~\ref{fig:orig_8_12mi_k973} presents an example of how the cycling
	path can change as the bicycle network is extended. The corresponding
	OD has a shortest path length of 4,250 meters (2.64 miles, $s_k$), and
	the maximum allowed length for a bicycle safe path is 6,375 meters
	(3.96 miles, $L_k=1.5s_k$).  On the current network, there does not
	exist a safe bicycle path that is sufficiently short.  However, the OD
	achieves a safe bicycle path of 6,009 meters (3.73 miles, $1.41 s_k$)
	when $B=12.8$ km of roads are improved.  Increasing the budget further
	to 19.2 km provides a shorter bicycle path of 4,976 meters (3.09
	miles, $1.17 s_k$).

	\subsection{The Benefit of Optimization}
	\label{sec:optimization solution advantage}
	
	This section compares the optimal improvement plans (OPT) to those
	obtained with the greedy heuristic (HEU). For instance, Figure
	\ref{fig:1.1 1.2 1.3 1.5 comparison - heu}, which includes 38.4 km (24 miles)
	HEU plans of four deviation factors, demonstrates that heuristic
	improvements focus less on constructing the backbone, compared to
	Figure \ref{fig:1.1 1.2 1.3 1.5 comparison}, and often unnecessarily
	provide multiple connections that serve a similar purpose.  The
	optimal and heuristic plans are compared on the percentage of
	potential cyclists and the average additional trip distance.  An
	example set of heuristic solutions for $R=1.1$ is presented by
	Figure~\ref{fig:examplemaps_HEU} in Appendix \ref{sec:appen_example maps}.
	
	\paragraph{Individual Travel Comparison}
	
	By definition, a heuristic plan cannot be better than the optimal plan, but the heuristic may improve some trips more than the optimal plan. For example, Table \ref{tab:opt vs heu} counts the number of travels with a smaller trip distance when $R=1.1$: it shows that, for all budgets, the optimal plans produces shorter routes more often than the heuristic. 
	
	\begin{spacing}{1}
		\begin{table}[!t]
			\begin{center}
				\begin{tabular}{ c c c c } 
					\toprule
					& \multicolumn{3}{c}{Travels Count}\\
					\cmidrule{2-4}
					B (km) & OPT map & HEU map & Equal\\ 
					\midrule
					6.4 & 89 & 43 & 907\\
					12.8 & 123 & 25 & 891\\
					19.2 & 142 & 7 & 890\\
					25.6 & 255 & 54 & 730\\
					32.0 & 321 & 38 & 680\\
					38.4 & 338 & 65 & 636\\
					44.8 & 275 & 69 & 695\\
					\bottomrule
				\end{tabular}
			\end{center}
			\caption{Travels with Smaller Penalty on Different Improvement Plans.} \label{tab:opt vs heu}
		\end{table}
		\doublespacing
	\end{spacing}
	
\begin{figure}[!t]
		\centering
		\includegraphics[width=0.55\textwidth]{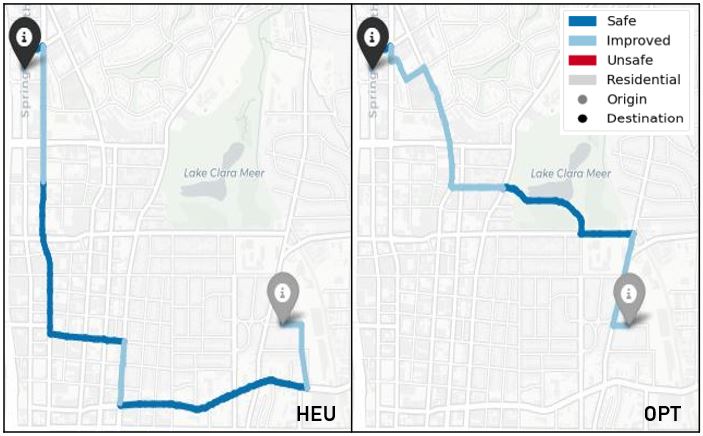}
		\caption{Example routes on OPT and HEU maps, $R=1.3$, $B=38.4$ km (24 mi).}
		\label{fig:different route}
	\end{figure}

	Figure~\ref{fig:different route} shows a single OD that is evaluated
	both on the heuristic plan and the optimal plan.  The corresponding OD
	has a shortest path length of 3,784 meters (2.35 miles, $s_k$), and
	the maximum allowed length for a bicycle safe path is $L_k=1.3s_k$.
	The optimal plan provides a short route of 3,973 meters (2.47 miles,
	$1.05s_k$) that is below the threshold distance.  The heuristic also
	produces a safe path; however, the path length is 5,286 meters (3.28
	miles, $1.39s_k$), which exceeds the threshold, so the
	heuristic plan still requires the outside option to serve the OD.

	\paragraph{Potential Cyclists}
	
	\begin{figure}[!t]
		\centering
		\includegraphics[width=0.6\textwidth]{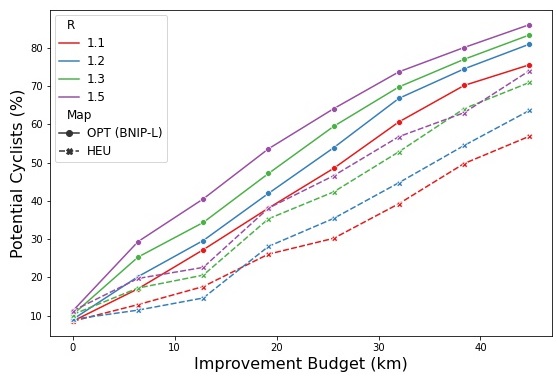}
		\caption{Percentage of Cyclists among Population.}
		\label{fig:pct_posutil}
	\end{figure}
	
	Figure~\ref{fig:pct_posutil} compares the heuristic plans and the optimal plans in terms of the number of potential cyclists. The results demonstrate that the optimal plans produce significant benefits in the number of potential cyclists. This is consistent over all budget values, and the difference may be more than 20\%. This is a compelling demonstration of the value of sophisticated optimization for infrastructure improvement.  
	
	\paragraph{Average Additional Trip Distance}
	
	\begin{figure}[!t]
		\centering
		\includegraphics[width=0.6\textwidth]{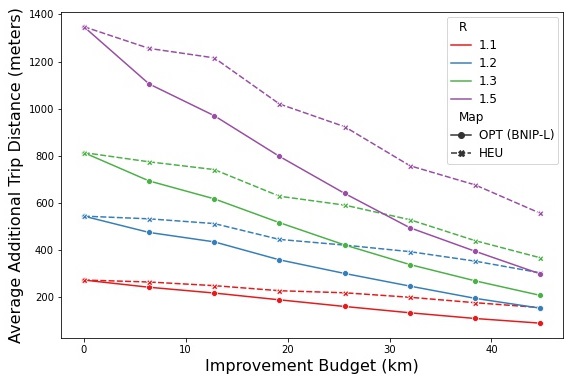}
		\caption{Average Additional Trip Distance.}
		\label{fig:utility_sum}
	\end{figure}
	
	Figure~\ref{fig:utility_sum} compares the average additional trip distance of the optimal and the heuristic plans for different values of $B$ and $R$. Again, the optimal plans produce significant benefits compared to the heuristics for all budget values. They parallel the improvements in potential cyclists and demonstrate the significant value of optimization for the BNIP. 
	
	\begin{figure}[!t]
		\centering
		\includegraphics[width=0.85\textwidth]{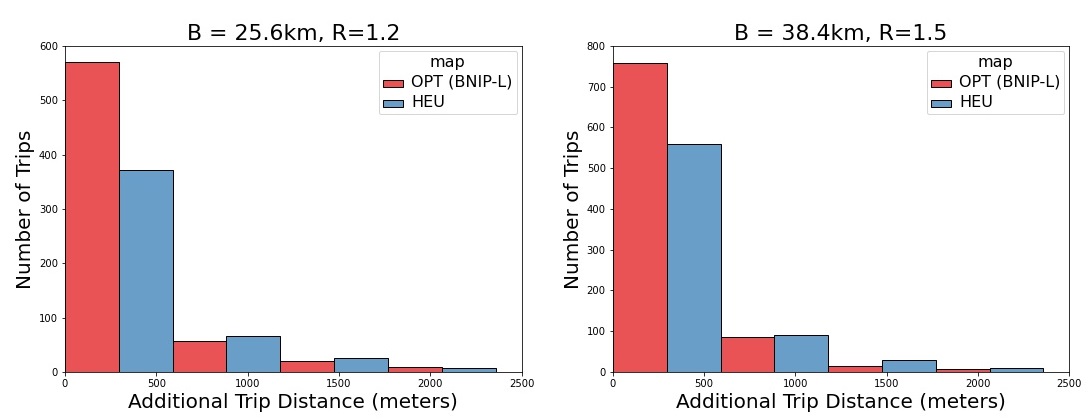}
		\caption{Distribution of Additional Trip Distances among Bicycle Travels.}
		\label{fig:util_dist}
	\end{figure}
	
	Figure \ref{fig:util_dist} shows the distribution of additional trip
	distances of safe routes for the optimal and heuristic plans under two
	example settings. Among the safe routes, the heuristic plans tend to
	produce fewer trips with smaller additional trip distances. In
	contrast, the optimization plans, which optimize both the number of
	potential cyclists and the average additional trip distances, have
	more trips with smaller additional trip distances as well as have more
	feasible bicycle travels.
	
	\subsection{Sequential Incremental Improvements}
	\label{sec:case study incremental improvement}
	
\begin{figure}[!t]
		\centering
		\includegraphics[width=\linewidth]{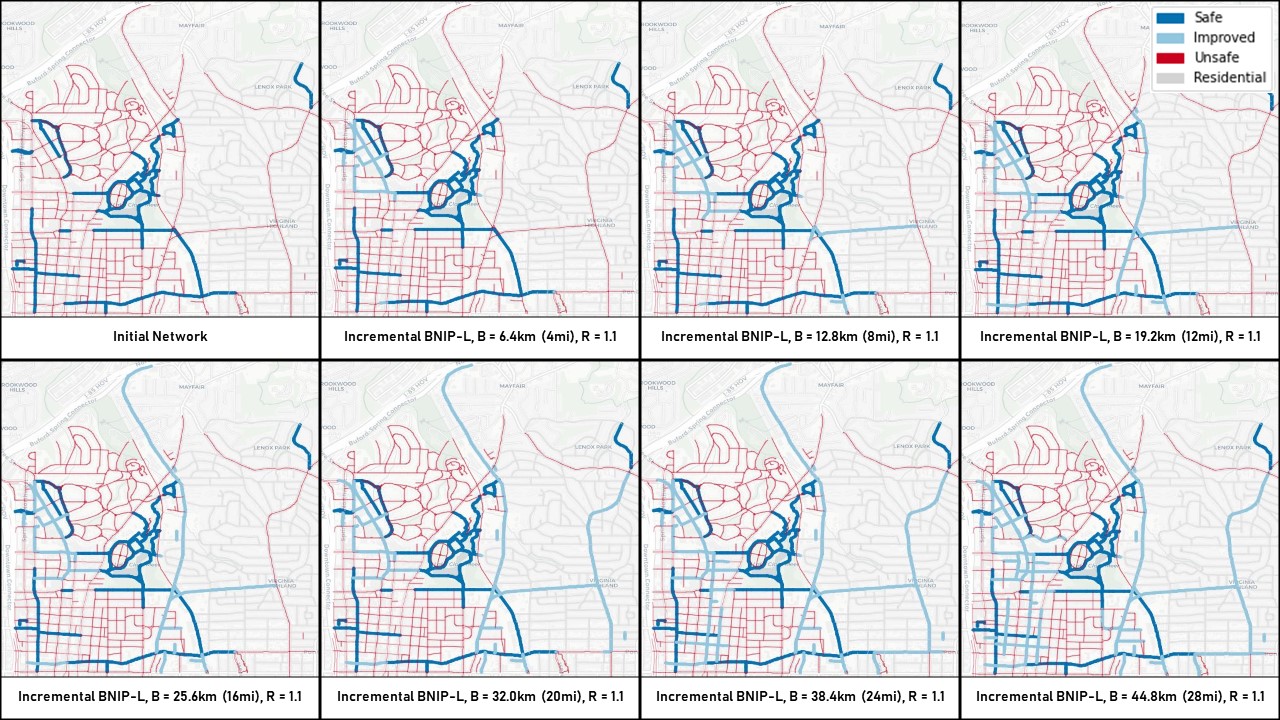}
		\caption{Sequential Incremental Improvement Plans using BNIP-L, $R=1.1$.}
		\label{fig:examplemaps_mono}
	\end{figure}

	In practice, policy makers cannot necessarily predict future
	availability of budget or other resources for the infrastructure
	development. Therefore, the improvement plans may be prepared over
	time and not in advance. Furthermore, there may be an incentive to
	implement solutions that are optimal in the short term, but may not
	necessarily be good in the long run. To investigate the effect of
	small myopic improvements, it is interesting to study the cumulative
	effect of a succession of individual improvement plans that optimally
	extend the bicycle infrastructure by 6.4 km at a time.
	Figure~\ref{fig:examplemaps_mono} presents one example of these plans,
	which illustrates how the sequential approach progressively produces
	the same overall structure as the optimal plan using $R=1.1$.
	
	\begin{figure}[!t]
		\centering
		\includegraphics[width=0.6\textwidth]{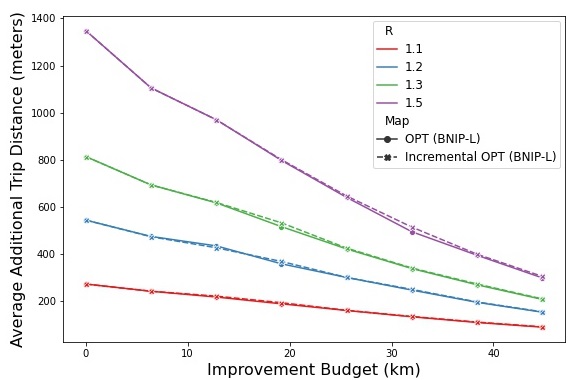}
		\caption{Average Additional Trip Distance for the Sequential Improvement Approach.}
		\label{fig:monotonicity}
	\end{figure}
	
	Figure~\ref{fig:monotonicity} shows that for all $R$, the additional average trip distances in the sequential approach is essentially similar to the overall optimal plans. This is of great practical importance, as it indicates that incrementally improving the network over time is practically identical to a strategic planning approach.
	
	\begin{spacing}{1}
		\begin{table}[!t]
			\begin{center}
				\begin{tabular}{ c c c c c c c c } 
					\toprule
					& \multicolumn{7}{c}{$B$}\\
					\cmidrule{2-8}
					Improvement Plans & 6.4 km & 12.8 km & 19.2 km & 25.6 km & 32.0 km & 38.4 km & 44.8 km \\
					\midrule
					\multirow{2}{*}{} 
					Strategic Planning & 40.49 & 65.07 & 93.90 & 122.12 & 148.88 & 173.10 & 192.99 \\
					Sequential Approach & 40.49 & 61.22 & 89.13 & 121.54 & 147.64 & 171.40 & 191.56 \\
					\midrule
					Difference & 0\% & 5.92\% & 5.07\% & 0.48\% & 0.83\% & 0.98\% & 0.74\% \\
					\bottomrule
				\end{tabular}
				\caption{Strategic Planning versus Incremental Improvements: Average Additional Trip Distance (m), $R=1.1$.} \label{tab:utility_score}
			\end{center}
		\end{table}
		\doublespacing
	\end{spacing}
	
	Table \ref{tab:utility_score} provides more details on the differences between the plans.
	The largest difference of additional trip distances between an optimal and a myopic plan is only 4950 meters in total, or 5 meters per travel on average.
	The optimal plans and the myopic plans do not lead to the same networks, but the difference in penalty is very small, especially as the network is improved further over time.
	
	\subsection{Uneven Improvement Costs}
	\label{sec:uneven improvement costs}
	Another practical consideration for the BNIP is to consider improvement costs that depend on the road properties, such as the number of lanes, pavement type, traffic volume, etc.
	Uneven improvement costs can be supported easily by replacing the distance parameters $d_{ij}$ in the budget constraint \eqref{eq:new:budget} by some other cost parameters $b_{ij}$ that combine multiple factors. Calculating the real cost for each way would take a significant amount of data, and is out of the scope of this paper. However, a sensitivity analysis based on the road type is provided in this section. 
	
	For this analysis, the ways $(i,j) \in W'$ are partitioned into three distinct sets based on the \emph{cycling-condition labels} retrieved from \citet{OpenStreetMap2020}.
	First, any road that has a cycling-related label, e.g., ``cycleway'' or ``bicycle'', but is known not to have a dedicated cycle lane is classified as a \emph{bike-friendly} road, and $b_{ij} = \frac{1}{2} d_{ij}$ is used as the weight.
	The remaining roads are ranked by the number of lanes and by how important the road is in the local system based on the following labels: motorway, trunk, primary, secondary, tertiary, residential, and unclassified.
	The roads with the first four labels are classified as \emph{significant} roads with higher improvement costs $b_{ij} = 2 d_{ij}$, and the others are set to incur distance values as costs ($b_{ij} = d_{ij}$).
	
\begin{figure}[!t]
		\centering
		\includegraphics[width=0.6\textwidth]{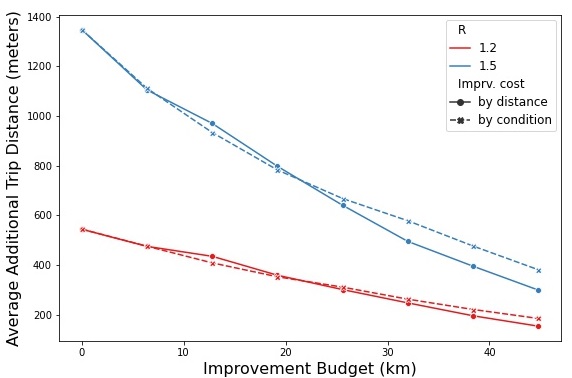}
		\caption{Average Additional Trip Distance with Uneven Improvement Costs.}
		\label{fig:uneven cost travel distance}
	\end{figure}

\begin{figure}[!t]
		\centering
		\includegraphics[width=0.6\textwidth]{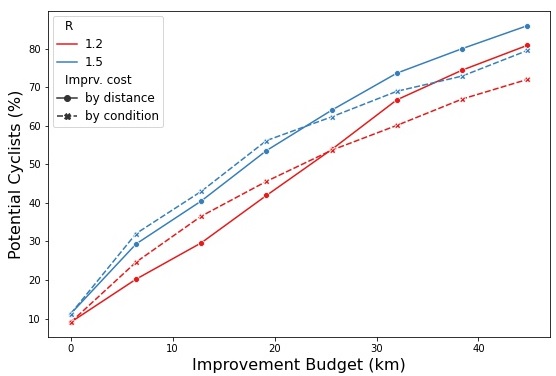}
		\caption{Percentage of Cyclists among Population with Uneven Improvement Costs.}
		\label{fig:uneven cost posutil}
	\end{figure}

	New improvement plans are prepared using these new costs for the same budget, under two deviation factors. 
	Figures~\ref{fig:examplemaps_bnip_unevencost 1.2} and~\ref{fig:examplemaps_bnip_unevencost 1.5} in Appendix~\ref{sec:appen_example maps} present such results.
	For the lower budget cases it is observed that the algorithm takes advantage of cheaper improvement costs by adding more lanes near the business area rather than improving major backbone roads.
	When the budget is increased, the results eventually converge to include all backbone roads, and then the maps are essentially identical to the original solutions.
	Figures \ref{fig:uneven cost travel distance} and \ref{fig:uneven cost posutil} further compare the additional average trip distances and the percentage of potential cyclists for the improvement plans designed by distance or by road condition.
	The trends are very similar, which suggests that the results are robust to uneven improvement costs.
	
	\section{Alternative Penalty Functions}
	\label{sec:comparing BNIP formulations}
	The results in the case study were generated with a linear penalty function, but other choices may lead to different bicycle network improvement plans.
	The BNIP supports different penalty functions, which makes it a flexible tool to use in practice.
	This section considers two alternative models aimed at minimizing the penalty of lost travelers, and maximizing the number of cyclists, respectively.
	The former results in a piecewise linear penalty function, and the latter requires a modification of the program formulation.
	Additional possibilities are discussed in Section~\ref{sec:discussion}.
	
	\subsection{Minimizing the Penalty of Lost Cyclists}
	This section presents a model that minimizes the penalty incurred by the travelers who choose \emph{not} to use the improved bicycle network.
	Compared to the linear penalty functions, the new model specifically focuses on potential cyclists that are lost to the system, rather than on the average additional trip distance.
	
	The model assumes that potential cyclists have a probability to drop out that increases with the deviation $u_k$ from the shortest safe trip.
	For short deviations $u_k \le 0.2s_k$, the probability is assumed to be zero, and no penalty is incurred.
	After that, the dropout probability increases linearly until the maximum trip length of $L_k = 1.5s_k$ is reached, at which point all travelers choose the outside option.
	The expected penalty for lost travelers at deviation $u_k \in [0.2 s_k, 0.5 s_k]$ follows from multiplying the dropout penalty of $u_k = L_k - s_k$ by the probability of dropping out.
	This results in the following penalty function:
	\begin{equation}
		f^P_k(u_k) =
		\begin{cases}
			0 & \textrm{ if } u_k \leq 0.2 s_k, \\
			\frac{5}{3} u_k - \frac{1}{3} s_k & \textrm{ if } 0.2 s_k \le u_k \leq 0.5 s_k.
		\end{cases}
	\end{equation}
	Note that $f^P_k$ is piecewise linear and convex. 
	The endpoint $f^P_k(0.5 s_k) = 0.5 s_k$ corresponds to a 100\% dropout rate and penalty $L_k - s_k = 0.5 s_k$.
	
	To use $f^P_k$ in the Benders decomposition algorithm, the master problem defines continuous $\upsilon$-variables, uses $f^P_k(\upsilon_k) = \upsilon_k$ in the objective, and includes the following two additional constraints:
	\begin{align}
		&\upsilon_k \geq 0 && \, \forall k \in T, \label{eq:pwlpenalty1}\\
		&\upsilon_k \geq \frac{5}{3} u_k - \frac{1}{3} s_k && \, \forall k \in T. \label{eq:pwlpenalty2}
	\end{align}
	The optimization with piecewise linear penalty function $f^P_k$ is labeled BNIP-P, and Figure~\ref{fig:examplemaps_bnip_pen} in Appendix~\ref{sec:appen_example maps} presents the full series of the BNIP-P improvement plans.
	
	\paragraph{Improvement Results}
	The BNIP-P improvement plans are compared with the BNIP-L plans for $R=1.5$, which have equal travel distance thresholds.
	\begin{spacing}{1}
		\begin{table}[!t]
			\centering
			\begin{tabular}{c c c c c c}
				\toprule
				& Network & \multicolumn{2}{c}{Average Additional} & \multicolumn{2}{c}{Cyclists}\\
				& Difference & \multicolumn{2}{c}{Distance (m)} & \multicolumn{2}{c}{Percentage}\\
				\cmidrule{3-6}
				B (km) & & BNIP-P & BNIP-L & BNIP-P & BNIP-L \\
				\midrule
				6.4 & 0\% & 1,105 & 1,105 & 29.36\% & 29.36\%\\
				12.8 & 4.20\% & 974 & 970 & 41.00\% & 40.52\%\\
				19.2 & 8.56\% & 972 & 797 & 53.90\% & 53.61\%\\
				25.6 & 8.66\% & 870 & 640 & 65.06\% & 64.10\%\\
				32.0 & 3.86\% & 515 & 494 & 75.65\% & 73.72\%\\ 
				38.4 & 3.00\% & 408 & 394 & 81.61\% & 80.08\%\\
				44.8 & 2.94\% & 672 & 298 & 87.68\% & 86.04\%\\
				\bottomrule
			\end{tabular}
			\caption{Network Improvement Comparison of BNIP-P and BNIP-L ($R=1.5$).} \label{tab:BNIP comparison network improv}
		\end{table}
		\doublespacing
	\end{spacing}
	Table~\ref{tab:BNIP comparison network improv} summarizes the improvements using each model.
	The two penalty functions produce very similar improvement plans; the percentage of network difference - a sum of percentages of unique improvement lengths on each improvement plan - shows that the two plans converge as more budgets are allowed.
	
	The deviation of an individual traveler is shorter for the BNIP-L,
	which is anticipated as the BNIP-P assigns equal penalty for short
	routes and focuses more on increasing the probability of bicycle
	participation.  However, the difference in cyclist percentages between
	the two optimized networks is very small.  This is an important
	observation, as it shows that the linear objective does not forfeit the advantage of
	the piecewise objective.
	
	\subsection{Maximizing the Number of Cyclists}
	
	Rather than assuming a probabilistic model for cyclists dropping out,
	it is also possible to directly maximize the number of potential
	cyclists, i.e., to maximize the number of ODs with a short and safe
	trip, or equivalently, minimize the use of the outside option $z_k$.
	The base formulation is modified to accommodate that objective, labeled as BNIP-Z, and the change is as follows:
	\begin{spacing}{1}
		\begin{mini!}
			%
			{}
			%
			{
				\sum_{k \in T} p_k z_k,
				\label{eq:new:objective-modified}
			}
			%
			{\label{formulation:new:BNIP-modified}}
			%
			{}
			%
			%
			\addConstraint
			{
				\sum_{(i,j)\in W'} d_{ij} y_{ij} \leq B, \label{eq:new:budget-modified}
			}
			{}
			{}
			\addConstraint
			{
				\sum_{(i,j)\in W} x^k_{ij} - \sum_{(j,i) \in W} x^k_{ji}
				=\begin{cases}
					1 - f_k & \textrm{ if } i = o_k\\
					f_k - 1 & \textrm{ if } i = d_k \\
					0 & \textrm{ otherwise}
				\end{cases} 
			}
			{}
			{\,\,\,\,\,\,\,\,\,\,\,
				\forall k\in T, i\in V,
				\label{eq:new:flow_balance-modified}
			}
			\addConstraint
			{x^k_{ij} \leq y_{ij}}
			{}
			{\forall k \in T, (i,j)\in W', \label{eq:new:upgraded_way_avail-modified}}
			\addConstraint
			{\sum_{(i,j) \in W} d_{ij} x^k_{ij} + (L_k + M_k) f_k \leq L_k + M_k z_k}
			{}
			{\forall k \in T, \label{eq:new:deviation-modified}}
			%
			%
			\addConstraint
			{y_{ij} \in \mathbb{B}}
			{}
			{\forall (i,j)\in W', \label{eq:new:s_binary-modified}}
			\addConstraint
			{x^k_{ij} \in \mathbb{B}}
			{}
			{\forall k \in T, (i,j)\in W, \label{eq:new:x_binary-modified}}
			\addConstraint
			{f_{k} \in \mathbb{B}}
			{}
			{\forall k \in T, \label{eq:new:y_binary-modified}}
			\addConstraint
			{z_{k} \in \mathbb{B}}
			{}
			{\forall k \in T. \label{eq:new:z_binary-modified}}
		\end{mini!}
		\doublespacing
	\end{spacing}
	
	Objective~\eqref{eq:new:objective-modified} minimizes the use of the outside option on the network.
	A binary variable $f_k$ is introduced for every trip $k \in T$ to ensure complete recourse, and the $f$-variables replace the $z$-variables in Constraints~\eqref{eq:new:flow_balance-modified}.
	That is, if no safe bicycle path exists, $f_k$ can be set to one instead.
	Constraints~\eqref{eq:new:deviation-modified} use the constant $M_k > 0$ to ensure that the outside-option variables $z_k$ are set correctly: preferably, the outside option is not used ($z_k=0$), which implies that $f_k=0$ and that the $x$-variables represent a safe and short bicycle path.
	If the outside option is used, setting $f_k=1$ removes this requirement.
	As a result, BNIP-Z indeed minimizes the use of the outside option.
	Finally, Constraints~\eqref{eq:new:y_binary-modified} define the newly added variables.
	
	It is notable that the Benders decomposition solution method can be used for the BNIP-Z.
	The modified master problem for the BNIP-Z is as follows:
	\begin{mini!}
		%
		{}
		%
		{
			\eqref{eq:new:objective-modified}, \notag
		}
		%
		{\label{formulation:new:master problem-modified}}
		{}
		%
		%
		\addConstraint
		{
			\eqref{eq:new:budget-modified}, \eqref{eq:new:s_binary-modified},  \eqref{eq:new:z_binary-modified}, \notag
		}
		{}
		{}
		\addConstraint
		{\Psi_k(y) \leq L_k + M_k z_k \quad \forall \, k \in T. \label{eq:new:master_cut-modified}\tag{\theequation}}
		{}
		{}
	\end{mini!}%
	The modified subproblem for the BNIP-Z is as follows:
	\begin{alignat}{4}
		\Psi_k(y) = \min \,&\,\,  \sum_{(i,j) \in W} d_{ij} x^k_{ij} + (L_k + M_k) f_k,  & & & & \\
		\text{s.t.} \,\,
		& \eqref{eq:new:flow_balance-modified}, \eqref{eq:new:upgraded_way_avail-modified},
		\eqref{eq:new:x_binary-modified}, \eqref{eq:new:y_binary-modified}. \,\,\, \notag &&
	\end{alignat}
	The new function calculating the shortest path, $\Psi_k(y)$, is introduced by the modified subproblem.
	The master problem now includes the $z$-variables for the minimization of Objective \eqref{eq:new:objective-modified}.
	The subproblem can be solved by linear programming as it is a standard minimum-cost flow problem.
	Pareto-optimal cuts and two-phase Benders can be used in the same way as before, and both $y$ and $z$-variables may be relaxed in phase one.
	
	Figure~\ref{fig:examplemaps_bnip_z} in Appendix~\ref{sec:appen_example maps} presents the full series of the BNIP-Z improvement plans using $L_k=1.5s_k$.
	
\begin{figure}[!t]
		\begin{adjustbox}{center}
			\includegraphics[width=\linewidth]{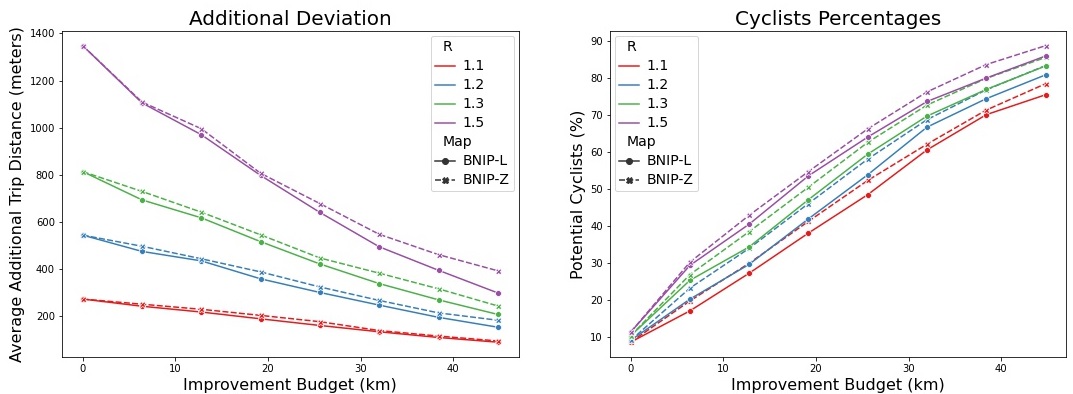}
		\end{adjustbox}
		\caption{Performances of the BNIP-L and BNIP-Z.}
		\label{fig:bnipz results}
	\end{figure}

	\paragraph{Improvement Results}
	Due to their objectives, the BNIP-L produces shorter deviations in average
	in its routes, and the BNIP-Z collects more potential cyclists.
	Nevertheless, Figure~\ref{fig:bnipz results} illustrates that the
	difference of performances between the two programs is
	inconsequential.  This shows that the use of the linear penalty, which
	examines both travel safety and proximity, does not sacrifice cyclists
	count to provide shorter deviation for riders.
	
	\subsection{Efficiency of Alternative Models}
\begin{figure}[!t]
		\begin{adjustbox}{center}
			\includegraphics[width=\linewidth]{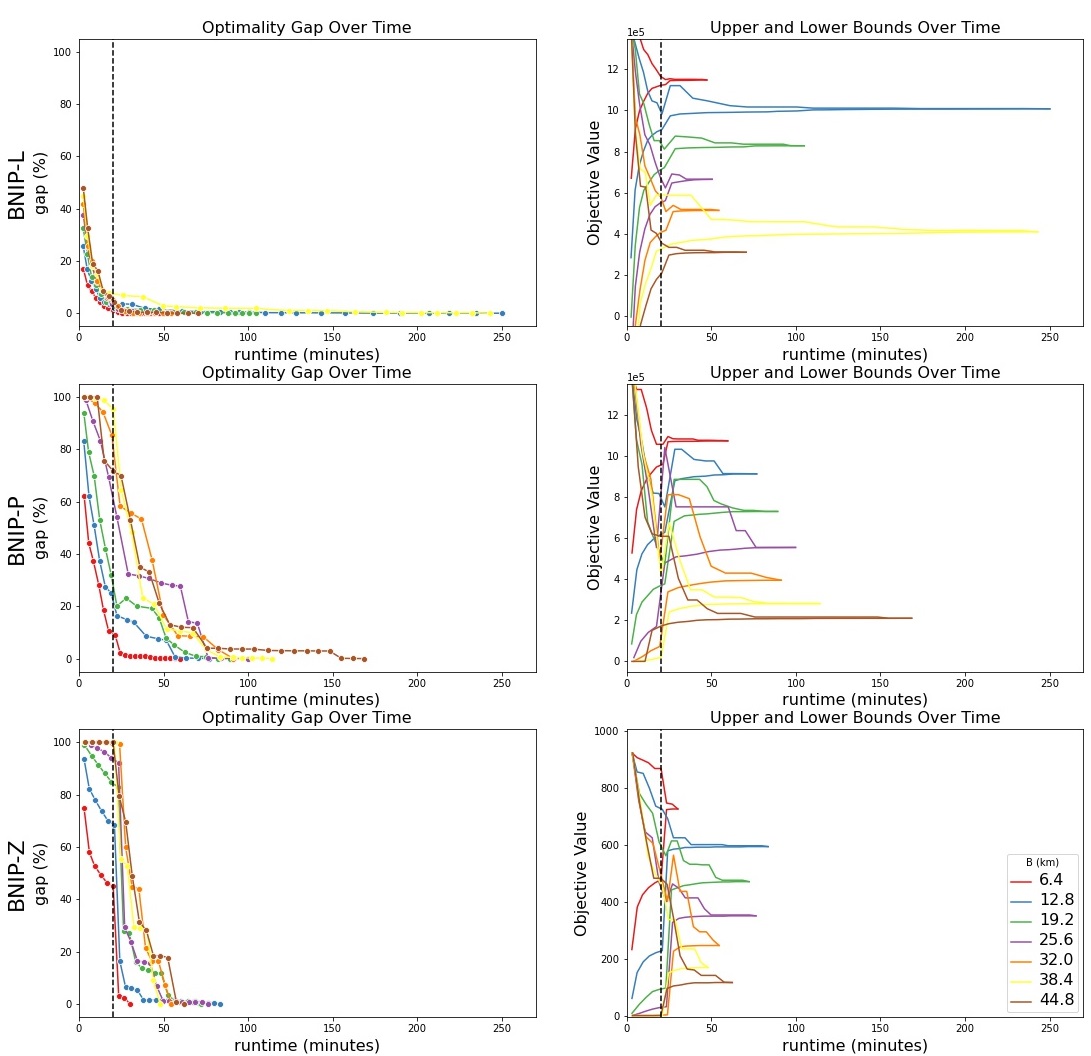}
		\end{adjustbox}
		\caption{Benders Decomposition Efficiency of the BNIP with Different Objectives, $L_k=1.5s_k$.}
		\label{fig:gap bnip comparison}
	\end{figure}

	To evaluate the methodological efficiency of the three BNIP
	formulations, the optimality gap and the upper and lower bounds of the
	Benders solutions are compared over time.  Using algorithm MW-McD,
	the BNIP-Z has the fastest computation times in many cases as shown in
	Figure \ref{fig:gap bnip comparison}. The BNIP-L also exhibits fast
	convergence to a small optimality gap for most cases.  This is
	noteworthy as BNIP-L considers both traveler safety and deviation.
	
	All three formulations reach optimality within a reasonable time, and
	they are shown to produce equally attractive improvement plans. Recall
	the BNIP-Z disregards travel distances, and the BNIP-P has a piecewise linear
	objective that is harder to interpret. The linear objective of the
	the BNIP-L reports the quality of the network in distance values, allowing
	for a more direct and meaningful evaluation of the network.
	
	\section{Discussion}\label{sec:discussion}
	The models in this paper focus on safety and travel distance, which are among the most important determinants for bicycle commuting \citep{HeinenEtal2010-CommutingByBicycle,CerveroEtAl2019-NetworkDesignBuilt,OspinaEtal2020-ColumbiaTravelDistance}.
	The case study is based on realistic data, and the conclusions are robust under uneven improvement costs (Section~\ref{sec:uneven improvement costs}) and alternative penalty functions (Section~\ref{sec:comparing BNIP formulations}).
	This section discusses how the models can be extended to accommodate
	additional data through alternative \emph{penalty functions}, \emph{cost structures}, and \emph{choice models}.
	
	\paragraph{Penalty Functions}
	The BNIP defines a term $p_k f_k(u_k)$ in the Objective~\eqref{eq:new:objective} for every trip $k \in T$.
	The parameter $p_k$ models the number of travelers completing this travel, but it can be used more generally to model any non-negative trip weights.
	This allows policy makers to assign more weight to certain areas or certain subsets of the population.
	From a model perspective, the only real requirements on the penalty function is that $f_k$ is non-decreasing (the requirement $f_k(0)=0$ is not restrictive).
	The penalty function is trip-specific, which means that external data (e.g., census data) can be used to help shape this function.
	
	From a computational perspective, the methods in this paper are expected to be effective for two general classes of penalty functions.
	The first class is that of convex penalty functions, which includes both BNIP-L (linear) and BNIP-P (piecewise linear).
	Linear and quadratic penalty functions can be handled directly by modern solvers (as for BNIP-L), and piecewise-linear functions can be modeled with additional variables and constraints (similar to \eqref{eq:pwlpenalty1}-\eqref{eq:pwlpenalty2} for BNIP-P).
	General convex functions can be handled with a classical cutting plane method \citep{KelleyJr.1960-CuttingPlaneMethod}.
	The second class is that of MIP representable functions, i.e., functions that can be modeled with mixed continuous and integer variables and additional constraints.
	BNIP-Z falls in this class, as it uses binary variables $z$ to indicate whether a traveler cycles or not.
	Modeling details are discussed by \citet{CroxtonEtAl2003-ComparisonMixedInteger}.
	It is worth noting that alternative penalty functions only affect how the master problem is solved, and the overall method stays the same.
	As such, similar computational performance is expected when the number of additional variables and constraints is small.
	
	\paragraph{Cost Structures}
	The BNIP includes two types of costs: the path length $l_k$ of trip $k \in T$, and the cost $d_{ij}$ for improving way $(i,j)\in W'$.
	Both are currently based on travel distance, but they can be used to reflect any combination of attributes.
	
	Path length can be replaced by a more general path cost.
	The \emph{path cost} $c_k$ is defined as the sum of a trip-specific constant $\alpha_k$ and trip-specific cost parameters $c^k_{ij} \ge 0$ for every way $(i,j) \in W$ on the path.
	The cost per way may combine any number of properties, including traffic stress, road gradient, activity density, mean rainfall, and even trip-specific socio-demographic attributes.
	Many of these attributes are discussed by \citet{CerveroEtAl2019-NetworkDesignBuilt}.
	The maximum acceptable length $L_k$ is replaced by the maximum acceptable cost $C_k$, accordingly.
	The formulation is updated by replacing
	\begin{equation}
		u_k \ge \sum_{(i,j) \in W} d_{ij} x^k_{ij} +  L_k z_k - s_k \quad \forall k \in T, \textrm{ by}\tag{\ref{eq:new:deviation}}
	\end{equation}
	\begin{equation}
		u_k \ge \sum_{(i,j) \in W} c^k_{ij} x^k_{ij} +  C_k z_k + \alpha_k \quad \forall k \in T.
	\end{equation}
	The original model reappears when the cost is distance ($c_{ij}^k = d_{ij}$, $C_k = L_k$) and the constant $\alpha_k =-s_k$ is used to calculate the deviation from the shortest path.
	Switching from path length to path cost makes the model more expressive and allows for including additional data, without affecting the solution method.
	
	A similar argument can be made for the cost of improving ways.
	The distance $d_{ij}$ in budget constraint~\eqref{eq:new:budget} may be replaced by a general budget cost, as is done in Section~\ref{sec:uneven improvement costs} to account for different road types.
	Furthermore, additional constraints may be added without affecting the solution method.
	For example, it can be enforced that the budget is spread out evenly over different areas.
	
	\paragraph{Choice Models}
	This paper uses a simple choice model based on distance: if path length $l_k \le L_k$, then trip $k\in T$ is completed by cycling, and if $l_k > L_k$, the outside option is used.
	Based on the discussion above, a more expressive choice model is also supported: if cost $c_k \le C_k$, then trip $k\in T$ is completed by cycling, and the outside option is used if $c_k > C_k$.
	
	More generally, the methods in this paper support \emph{any} classifier that predicts cycling when $c_k \le C_k$ for some cut-point $C_k$, and no cycling otherwise.
	This includes logit and machine learning models \citep{ZhaoEtAl2019-ModelingHeterogeneityMode}.
	For example, the logistic regression model \citep{KleinbaumKlein2010-LogisticRegression} is given by
	\begin{equation}
		\pi(c_k) = \frac{1}{1 + e^{c_k}},
	\end{equation}
	where $\pi(c_k)$ is the probability of using the outside option.
	Fitting the model amounts to combining the different attributes into a cost $c_k$ that best explains traveler behavior.
	The model predicts cycling when the probability of using the outside option is low, i.e., $\pi(c_k) \le \Pi$ for some value $\Pi$.
	This probability cut-point can be translated into a cut-point for the traveler cost:
	\begin{equation}
		\pi(c_k) \le \Pi \Longleftrightarrow c_k \le \log\left(\frac{\Pi}{1-\Pi}\right) = C_k.
	\end{equation}
	This results in a setting that is directly supported by the methods in this paper.

	\section{Conclusion}
	\label{sec:conclusion}
	
	Cycling brings many benefits to both the cyclists and to society as a
	whole, and the emergence of e-bikes may make this transportation mode
	attractive for a larger population segment. However, safety is a
	critical issue faced by commuters when deciding their transportation mode. 
	This paper considered the problem of improving the bicycle infrastructure to allow more people to travel by bicycle.
	This optimization problem was formalized as the Bicycle Network Improvement Problem (BNIP).
	As opposed to the literature, the BNIP supports a budget for improvement, provides completely safe routes, allows full flexibility in routing cyclists, and the solution approach is exact.
	Solving the BNIP directly is
	computationally intractable for large instances, so the paper
	presented a Benders decomposition to remedy this issue
	by exploiting the problem structure and considering each rider
	independently in the subproblems. 
	
	The paper demonstrated the effectiveness of the method on an in-depth
	case study for Midtown Atlanta, based on real transportation data of
	white-collar commuters from Virginia-Highland.  The computational
	results show that the proposed Benders decomposition algorithm with
	Pareto-optimal cuts and two-phase Benders is very effective in solving the realistic case
	study instances.  Further analysis revealed that the optimal bicycle
	network improvement plans for Midtown Atlanta are very powerful in
	providing access to safe and short bicycle routes. The increase in
	the number of travelers with access to a safe and short trip is almost linear in the available budget,
	indicating that more investments in bicycle infrastructure may keep
	attracting additional commuters to switch to cycling.  The Benders
	decomposition method was compared to a greedy heuristic, and shown to lead to
	significantly better plans, which shows the value of optimization to produce mathematically optimal solutions.
	
	The paper also considered practical aspects of the bicycle network extension.
	It showed that repeated myopic
	extensions of the network leads to an almost optimal result in the
	long run.  This is of great practical importance, as it indicates that
	myopically improving the network over time is practically identical to
	in-advance planning. 
	It was shown that changing the road improvement costs to take the type of road into account did not affect the main conclusions of this paper.
	In addition, the paper demonstrated that the results
	are robust with respect to different objective functions.
	
	Future work may incorporate additional data and more complicated choice models into the BNIP.
	Section~\ref{sec:discussion} provides guidance on how this may be done.
	Another interesting direction is to investigate how to extend bicycle network improvement
	to multi-modal
	transit systems, where the goal is not necessary to offer a bicycle
	path for the complete trip, but to let bicycles play a role in
	addressing the ubiquitous first/last mile problem.
	
	
	\section{Acknowledgements}\label{sec:acknowledgements}
	This research is partly supported by NSF Leap HI proposal NSF-1854684.
	
	\clearpage
	\bibliographystyle{plainernat}
	\bibliography{LimEtAl2022}
	
	
	\clearpage
	\pagestyle{plain}
	\pagenumbering{Roman}
	\appendix
	
	\section{Bicycle Network Improvement Plans} \label{sec:appen_example maps}
\begin{sidewaysfigure}
		\centering
		\includegraphics[width=\linewidth]{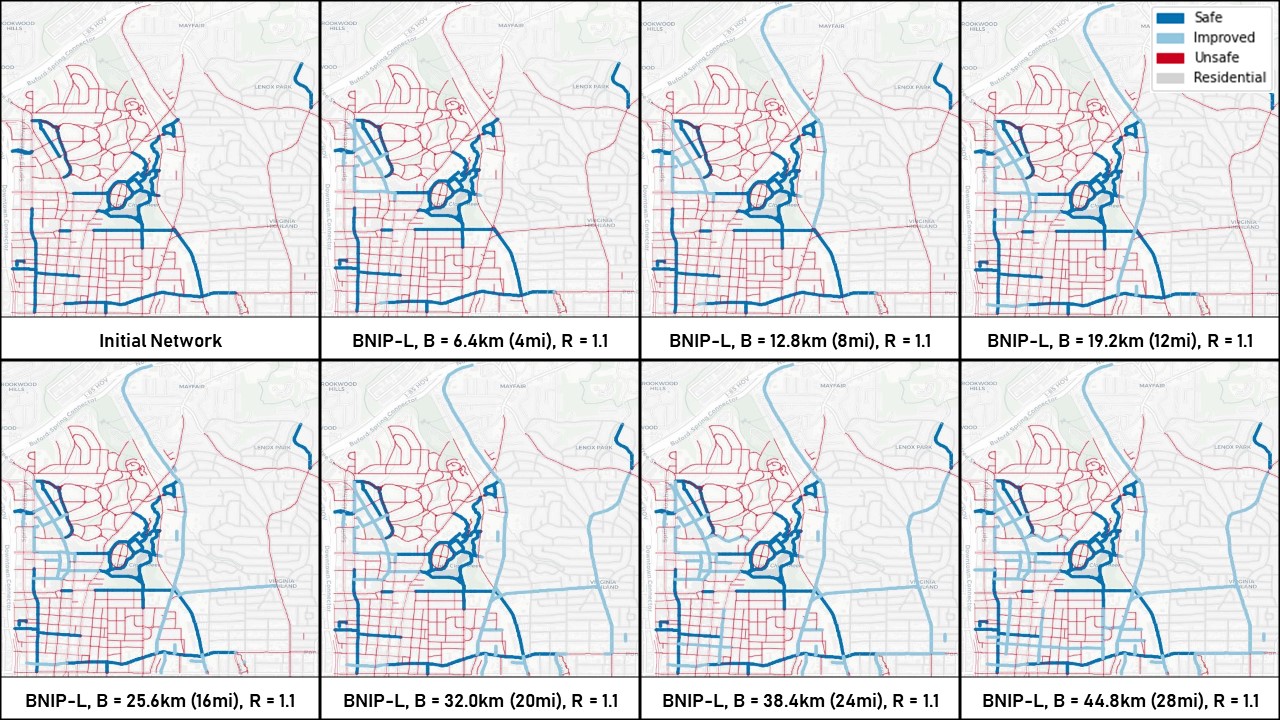}
		\caption{Optimal Bicycle Network Improvement Plans using BNIP-L, $R=1.1$.}
		\label{fig:examplemaps}
	\end{sidewaysfigure}

\begin{sidewaysfigure}
		\centering
		\includegraphics[width=\linewidth]{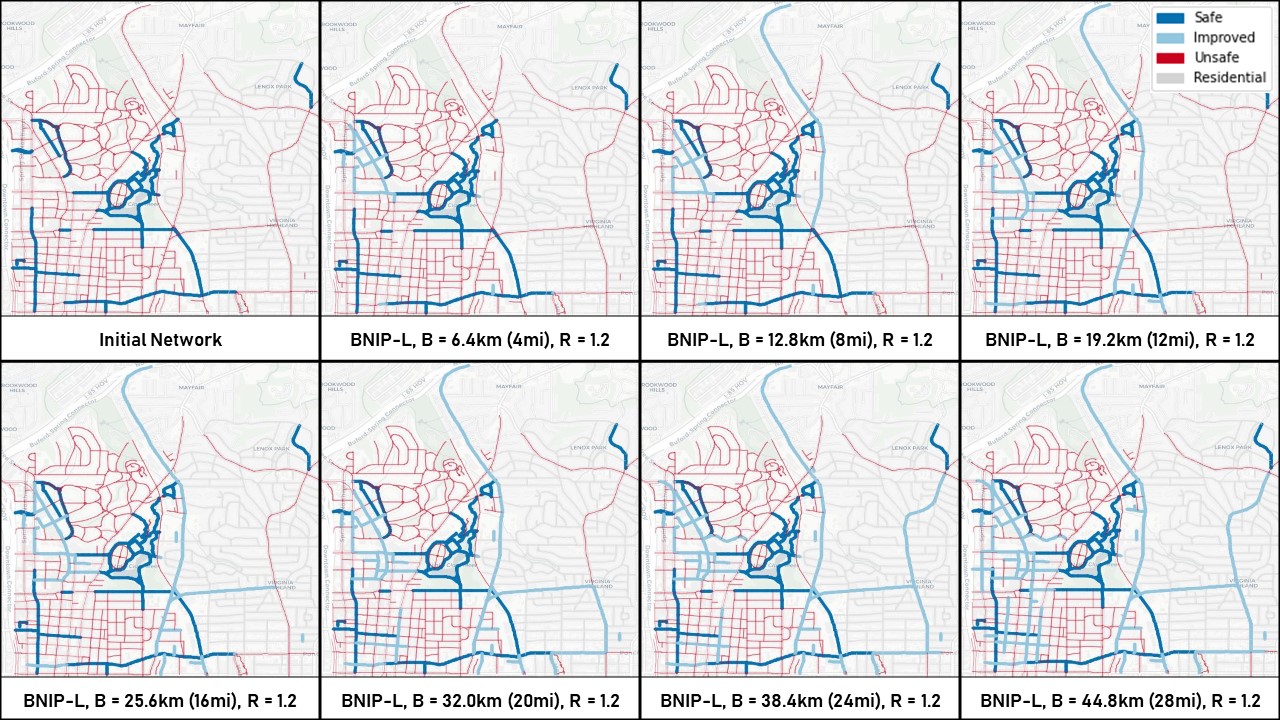}
		\caption{Optimal Bicycle Network Improvement Plans using BNIP-L, $R=1.2$.}
		\label{fig:examplemaps_1.2}
	\end{sidewaysfigure}

\begin{sidewaysfigure}
		\centering
		\includegraphics[width=\linewidth]{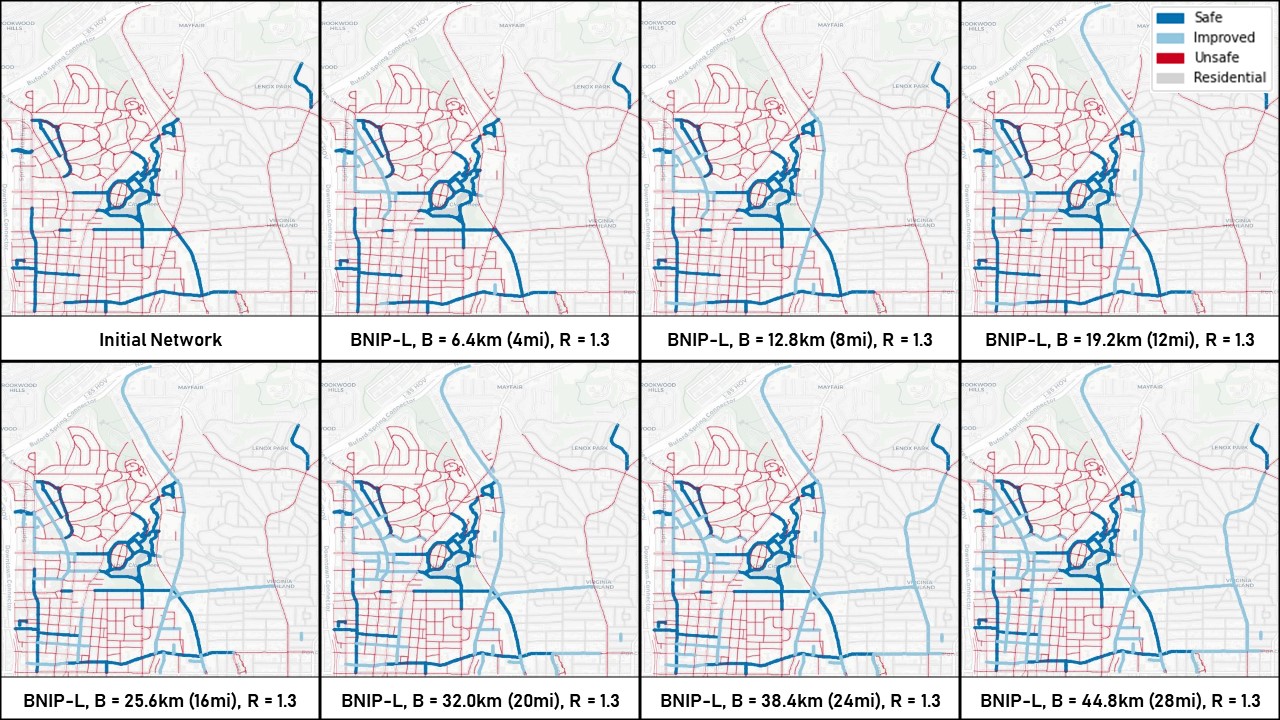}
		\caption{Optimal Bicycle Network Improvement Plans using BNIP-L, $R=1.3$.}
		\label{fig:examplemaps_1.3}
	\end{sidewaysfigure}

\begin{sidewaysfigure}	
		\centering
		\includegraphics[width=\linewidth]{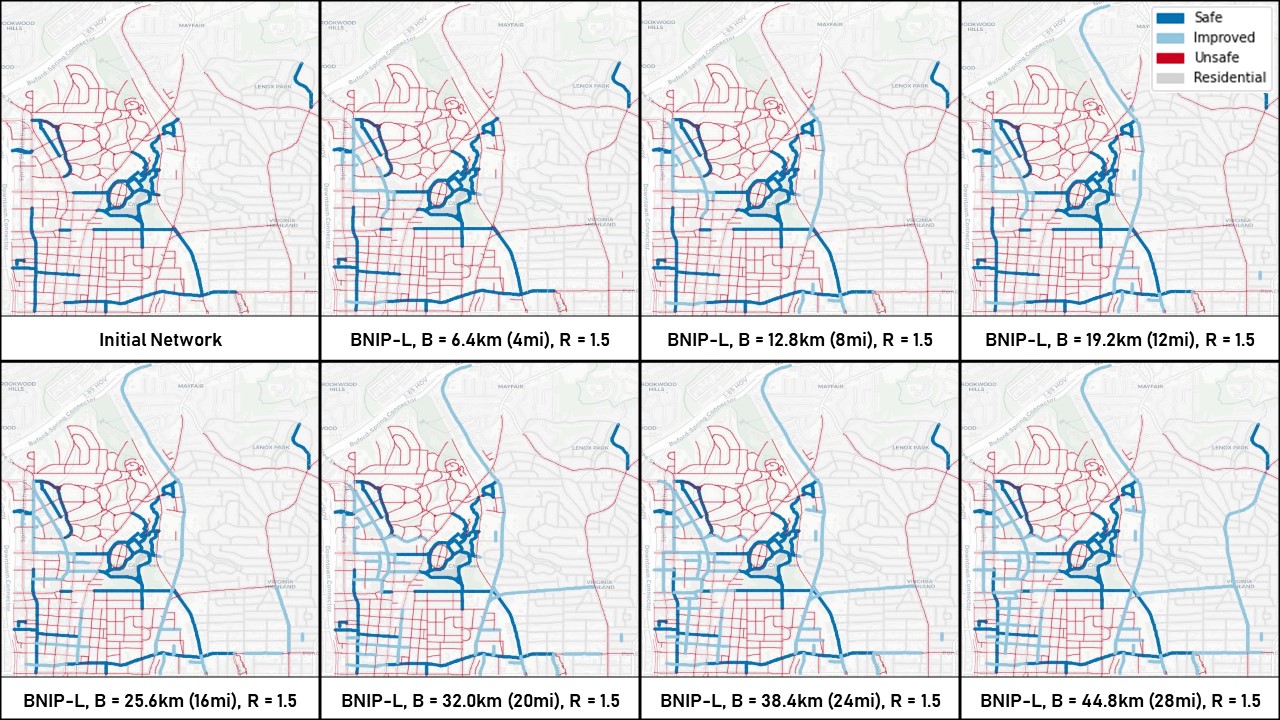}
		\caption{Optimal Bicycle Network Improvement Plans using BNIP-L, $R=1.5$.}
		\label{fig:examplemaps_1.5}
	\end{sidewaysfigure}

\begin{sidewaysfigure}
		\centering
		\includegraphics[width=\linewidth]{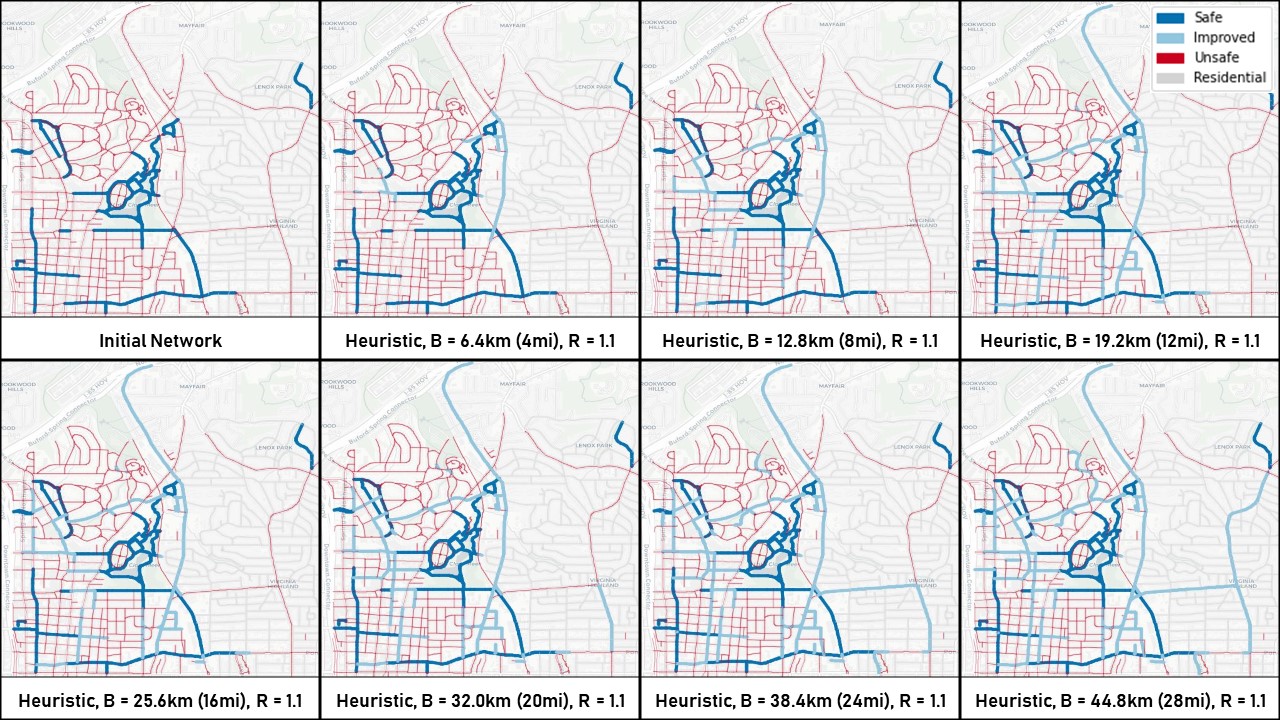}
		\caption{Heuristic Bicycle Network Improvement Plans, $L_k=1.1s_k$.}
		\label{fig:examplemaps_HEU}
	\end{sidewaysfigure}

\begin{sidewaysfigure}
		\centering
		\includegraphics[width=\linewidth]{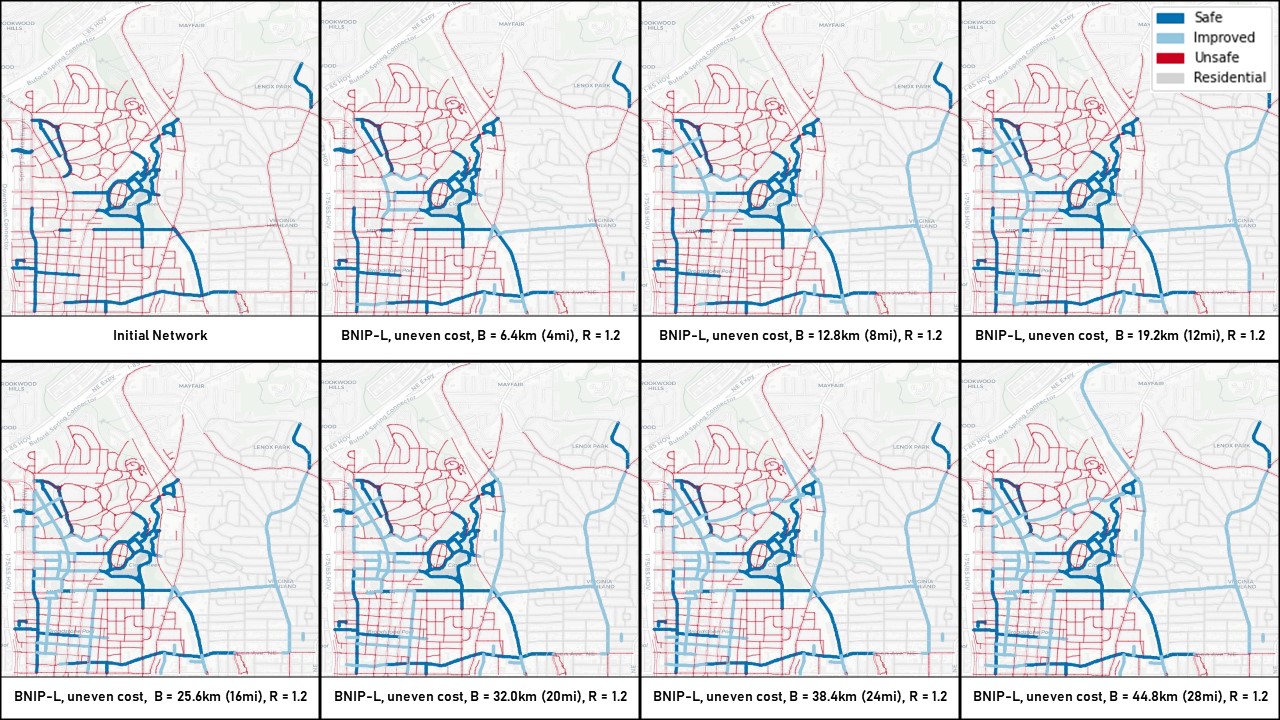}
		\caption{Optimal Bicycle Network Improvement Plans using BNIP-L and Uneven Improvement Costs, $R=1.2$.}
		\label{fig:examplemaps_bnip_unevencost 1.2}
	\end{sidewaysfigure}

\begin{sidewaysfigure}
		\centering
		\includegraphics[width=\linewidth]{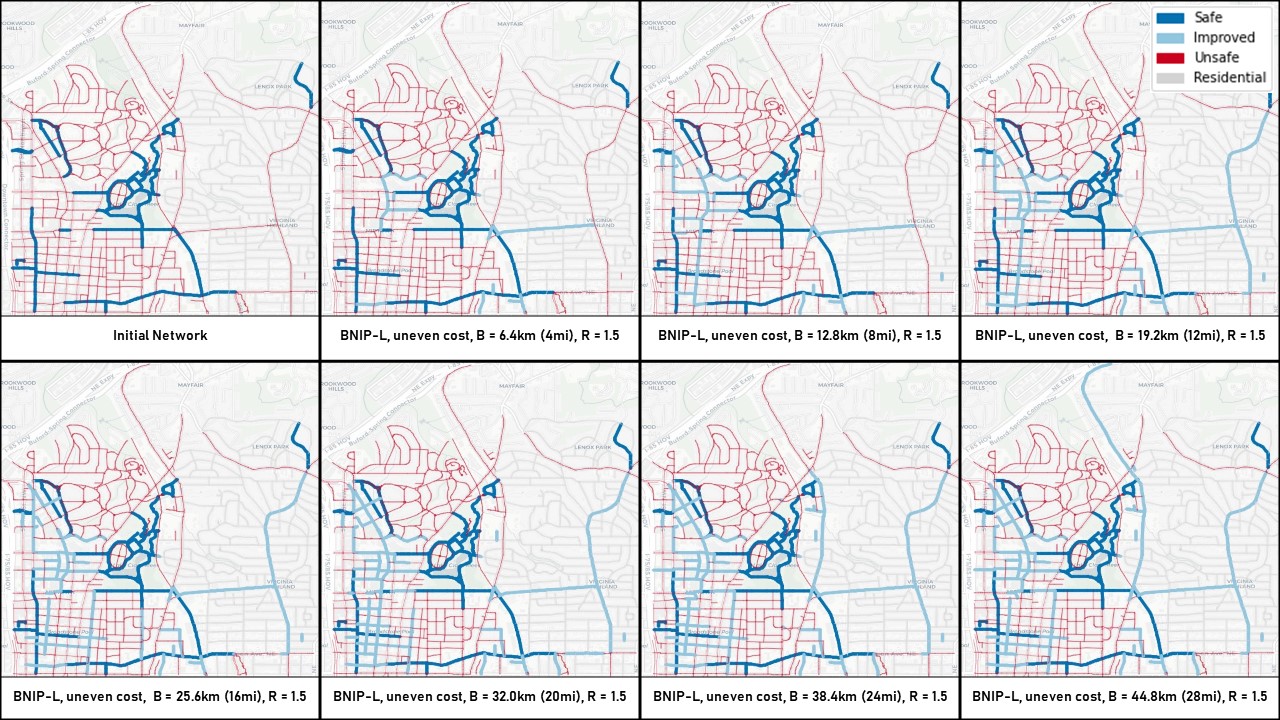}
		\caption{Optimal Bicycle Network Improvement Plans using BNIP-L and Uneven Improvement Costs, $R=1.5$.}
		\label{fig:examplemaps_bnip_unevencost 1.5}
	\end{sidewaysfigure}

\begin{sidewaysfigure}
		\centering
		\includegraphics[width=\linewidth]{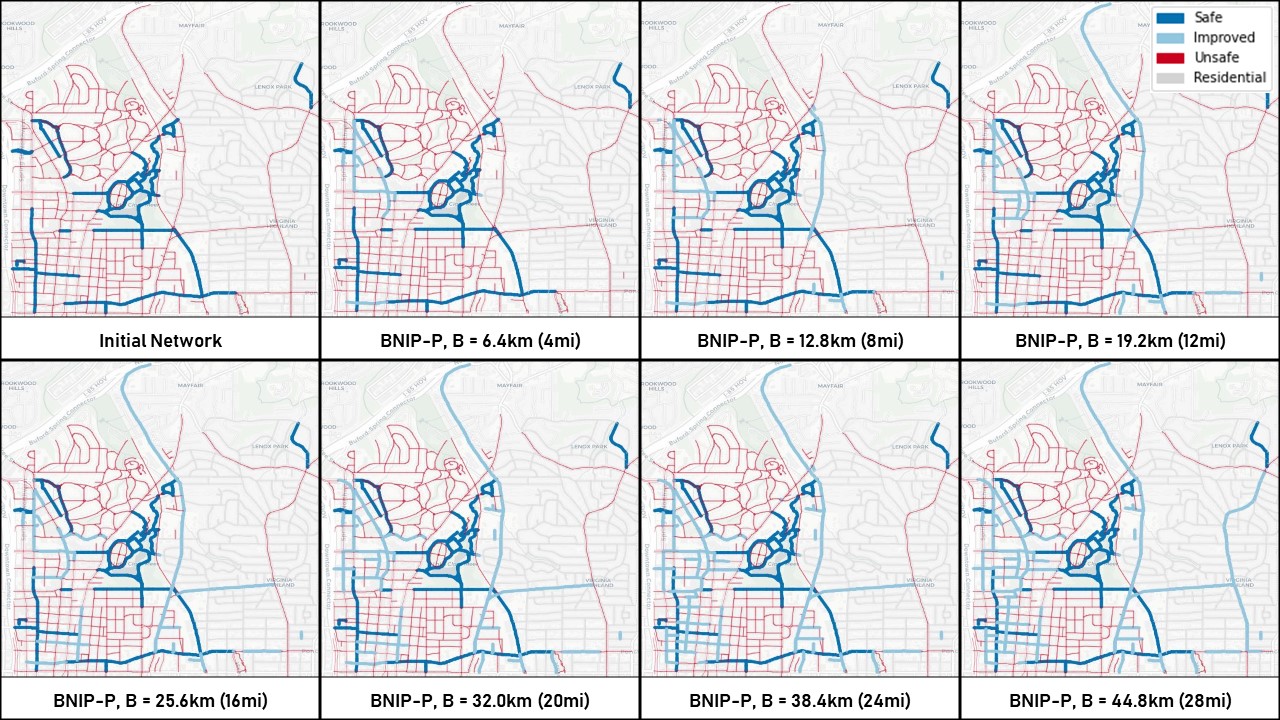}
		\caption{Optimal Bicycle Network Improvement Plans using BNIP-P.}
		\label{fig:examplemaps_bnip_pen}
	\end{sidewaysfigure}

\begin{sidewaysfigure}
		\centering
		\includegraphics[width=\linewidth]{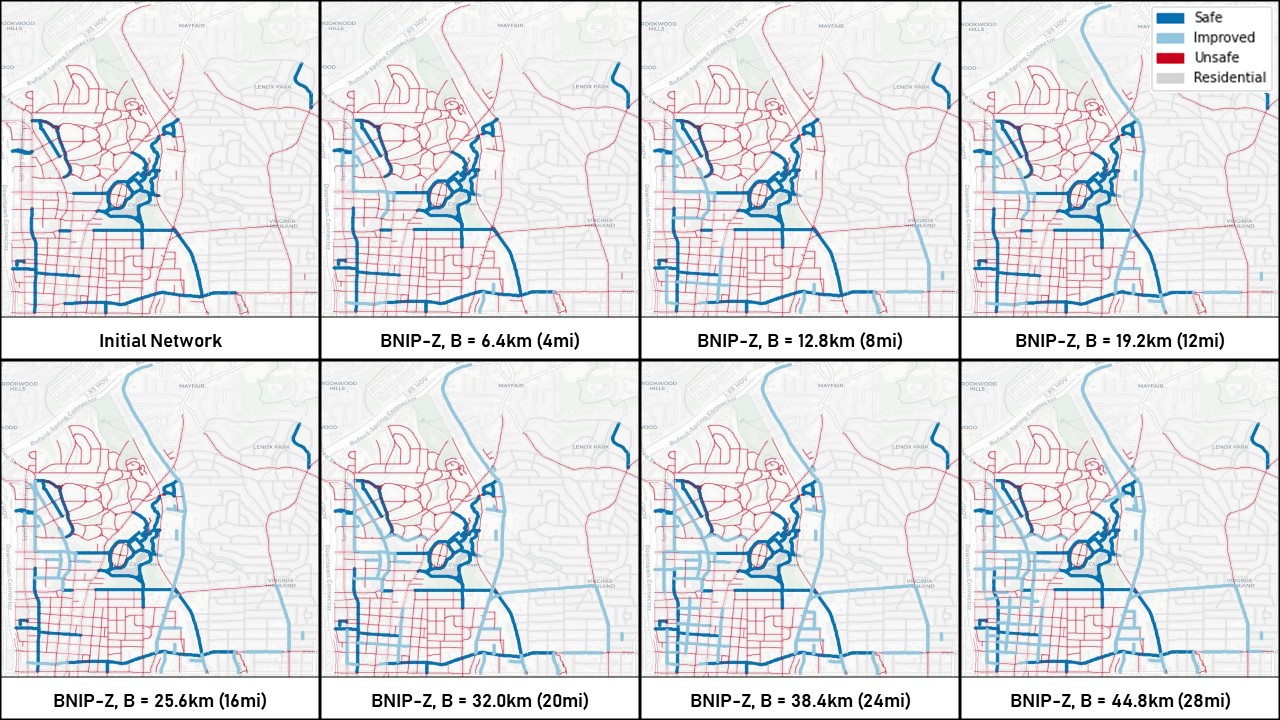}
		\caption{Optimal Bicycle Network Improvement Plans using BNIP-Z, $L_k=1.5s_k$.}
		\label{fig:examplemaps_bnip_z}
	\end{sidewaysfigure}

	This section presents all bicycle network improvement plans created for the case study.
	In Figures \ref{fig:examplemaps}-\ref{fig:examplemaps_bnip_z}, existing bicycle infrastructure is indicated by blue roads, and light-blue roads indicate the proposed expansion.
	Red roads remain unsafe, and gray roads are residential roads that are safe to use without improvement.
	
\end{document}